\documentclass[twoside,11pt]{article}

\usepackage{jmlr2e,times}
\usepackage[english]{babel}
\usepackage[latin1]{inputenc}
\usepackage{bbm,amssymb,amsmath,amsfonts,graphicx,fullpage,ifthen,twoopt,algorithm,algorithmic,color}
\def\Rset{\mathbb{R}}

\def\E{\mathbb{E}}
\def\P{\mathbb{P}}
\def\T{\mathcal{T}}

\def\rme{\operatorname{e}}
\newcommand{\1}[1]{\mathbbm{1}_{\left\{#1\right\}}}
\DeclareMathOperator*{\argmax}{arg\,max}

\def\ps{P}\def\Ps{\P_\pi}\def\mus{\mu}\def\Es{\E_\pi}

\def\F{\mathcal{F}}\def\A{\mathcal{A}}
\def\eqdef{\triangleq}

\newcommand{\Reward}[2]{S_{#1}(#2)}
\newcommand{\EReward}[2]{M_{#1}(#2)}
\newcommand{\MReward}[3][{}]{\ifthenelse{\equal{#1}{}}{\bar{X}_{#2}(#3)}{\bar{X}_{#2}(#3, #1)}}
\newcommandtwoopt{\Nplay}[4][{}][{}]{
	\ifthenelse{\equal{#1}{}}{
		\ifthenelse{\equal{#2}{}}
		{N_{#3}(#4)}
		{N_{#3}(#4, #2)}
	}{
		\ifthenelse{\equal{#1}{sum}}{
			\ifthenelse{\equal{#2}{}}
			{n_{#3}(#4)}
			{n_{#3}(#4, #2)}
		}{
			\ifthenelse{\equal{#2}{}}
				{N^{#1}_{#3}(#4)}
				{N^{#1}_{#3}(#4, #2)}
		}
	}
}
\newcommand{\Nbadplayed}[2]{\tilde{N}_{#2}(#1)}
\newcommand{\Ic}[3]{c_{#2}(#3, #1)}
\newcommand{\DUCB}[1]{#1}
\newcommand{\SWUCB}[1]{#1}
\newcommand{\CPE}[2]{\E\left[ \left. #1 \right| #2 \right]}
\newcommandtwoopt{\Np}[3][{}][{}]{
\ifthenelse{\equal{#1}{}}{
  \ifthenelse{\equal{#2}{}}
    {N(#3)}
    {N^{#2}(#3)}
}{
  \ifthenelse{\equal{#2}{}}
    {N_{#1}(#3)}
    {N_{#1}^{#2}(#3)}}
}
\def\eqsp{\;}
\def\KL{\mathrm{D}}
\def\T{\mathcal{T}}

\begin{document}
\title{On Upper-Confidence Bound Policies for Non-Stationary Bandit Problems}

\author{\name Aur\'elien Garivier \email aurelien.garivier@telecom-paristech.fr \\
\addr  Institut Telecom / Telecom ParisTech / Laboratoire LTCI / CNRS UMR 5141 \\
46 rue Barrault, 75634 Paris Cedex 13
\AND
\name Eric Moulines \email eric.moulines@telecom-paristech.fr \\
\addr  Institut Telecom / Telecom ParisTech / Laboratoire LTCI / CNRS UMR 5141 \\
46 rue Barrault, 75634 Paris Cedex 13}

\editor{unknown}

\maketitle

\begin{abstract}
Multi-armed bandit problems are considered as a paradigm  of the  trade-off between exploring the environment to find profitable actions and exploiting what is already known. In the stationary case, the distributions of the rewards do not change in time, \emph{Upper-Confidence Bound} (UCB) policies, proposed in \cite{agrawal:1995} and later analyzed in \cite{auer02finitetime}, have been shown to be rate optimal.

A challenging variant of the MABP is the non-stationary bandit problem where the gambler must decide which arm to play while facing the possibility of a changing environment.
In this paper, we consider the situation where the distributions of rewards remain constant over epochs and change at unknown time instants. We analyze two algorithms:  the \emph{discounted UCB} and the \emph{sliding-window UCB}. We establish for these two algorithms an upper-bound for the expected regret by upper-bounding the expectation of the number of times a suboptimal arm is played.
For that purpose, we derive a Hoeffding type inequality for self normalized deviations with a random number of summands. We establish a lower-bound for the regret in presence of abrupt changes in  the arms reward distributions. We show that the discounted UCB and the sliding-window UCB both match the lower-bound up to a logarithmic factor.
\end{abstract}

\begin{keywords}
Multi-armed bandit, reinforcement learning, deviation inequalities, non-stationary environment
\end{keywords}

\section{Introduction}
Multi-armed bandit (MAB) problems, modelling allocation issues under uncertainty, are fundamental to stochastic decision theory.
The archetypal MAB problem may be stated as follows: there is a bandit with $K$ independent arms. At each time step, the player can play only one arm and receive a reward. In the stationary case, the distribution of the rewards are initially unknown, but are assumed to remain constant
 during all games.  The player iteratively plays one action (pulls an arm) per round, observes the associated reward, and decides on the action for the next iteration. The goal of a MAB algorithm is to minimize  the expected regret over $T$ rounds, which is defined as the expectation of the difference between the total reward obtained by playing the best arm and the total reward obtained by using the algorithm (or \emph{policy}).
 The minimization of the regret is achieved by balancing \emph{exploitation},  the use of acquired information, with \emph{exploration}, acquiring new information. If the player always plays the arm which he currently believes to be the best, he might miss to identify another arm having an actually higher expected reward. On the other hand, if the gambler explores too often the environment to find profitable actions, he will fail to accumulate as many rewards as he could.
 For several algorithms in the literature (e.g.~\cite{lai:robbins:1985,agrawal:1995}), as the number of plays $T$ goes to infinity, the expected total reward asymptotically approaches that of playing a policy with the highest expected reward, and the regret grows as the logarithm of $T$. More recently, finite-time bounds for the regret have been derived (see \cite{auer02finitetime,audibert:munos:szepesvari:2007}).

Though the stationary formulation of the MABP allows to address exploration versus exploitation challenges in a intuitive and elegant way, it may fail to be adequate to model an evolving environment where the reward distributions undergo changes in time.
As an example, in the cognitive medium radio access problem \cite{lai:2007}, a user  wishes to opportunistically exploit the availability of an empty channel in a multiple channels system; the reward is the availability of the channel, whose distribution is unknown to the user. Another application is real-time optimization of websites by targetting relevant content at individuals, and maximize the general interest by learning and serving the most popular content (such situations have been considered in the recent Exploration versus Exploitation (EvE) PASCAL challenge by \cite{hartland:gelly:baskiotis:teytaud:sebag:2006}, see also \cite{koulouriotis:xanthopoulos:2008} and the references therein).
These examples illustrate the limitations of the stationary MAB models. The probability that a  given channel is available  is likely to change in time. The news stories a visitor of a website is most likely to be interested in vary in time.

To model such situations, we need to consider non-stationary MAB problems, where distributions of rewards may change in time. We show in the following that, as expected, policies tailored for the stationary case fail to track changes of the best arm.
In this paper, we consider a particular non-stationary case where the distributions of the rewards undergo abrupt changes. We derive a lower-bound for the regret of any policy, and we analyze two algorithms: the Discounted UCB (Upper Confidence Bound) proposed by Koczis and Szepesv\'{a}ri and the Sliding Window UCB we introduce. We show that they are almost rate-optimal, as their regret almost matches a lower-bound.

\subsection{The stationary MAB problem}
At each time $s$, the player chooses an arm $I_s\in\{1,\dots,K\}$ to play according to a (deterministic or random) policy $\pi$ based on the sequence of past plays and rewards, and obtains a reward $X_s(I_s)$\footnote{Note that we use here the convention that the reward after at time $s$ if the $i$-th arm is played is supposed to be $X_s(i)$ and not the $N_s(i)$-th reward in the sequence of rewards for arm $i$, where $N_s(i)$ denotes the number of time the arm $i$ has been played up to time $s$; while this convention makes no difference in the stationary case, because the distribution of the rewards are independent, it is meaningful in the non-stationary case, since the distribution of the arm \emph{may change} even if the arm has not been played. These models can be seen as a special instance of the so-called \emph{restless} bandit, proposed by \cite{whittle:1988}.}.
The rewards $\{X_s(i) \}_{s \geq 1}$ for each arm $i\in\{1,\dots,K\}$ are modeled by a sequence of independent and indentically distributed (i.i.d.) random variables from a distribution unknown to the player. We denote by $\mu(i)$ the expectation of the reward $X_1(i)$.

The optimal (oracle) policy $\pi^*$ consists in always playing the arm $i^*\in \{1,\dots,K\}$ with largest expected reward $$\mu(*)=\max_{1\leq i\leq K}\mu(i) \eqsp, \quad i^* = \argmax_{1\leq i\leq K}\mu(i) \eqsp.$$
The performance of a policy $\pi$ is measured in terms of \emph{regret} in the first $T$ plays,
which is defined as the expected difference between the total rewards collected by the optimal policy $\pi^*$ (playing at each time instant the arm $i^*$ with the highest expected reward) and the total rewards collected by the policy $\pi$.

Denote by $\Np[t][]{i}=\sum_{s=1}^t\1{I_s=i}$ the number of times arm $i$ has been played in the $t$ first games.
The expected regret after $T$ plays may be expressed as:
\begin{equation*}
\E_\pi\left[\sum_{t=1}^T \left\{\mu(*)-\mu(I_t)\right\} \right]
=\sum_{i\neq i^*} \left\{ \mu(*)-\mu(i)\right\} \E_\pi\left[\Np[T][]{i}\right],
\end{equation*}
where $\E_\pi$ the expectation under policy $\pi$.

Obviously, bounding the expected regret after $T$ plays essentially amounts to controlling the expected number of times a sub-optimal arm is played.
In their seminal paper, \cite{lai:robbins:1985} consider stationary MAB problem, in which the distribution of rewards was taken from a one-dimensional parametric family (each being associated with a different value of the parameter, unknown to the player). They have proposed a policy achieving a logarithmic regret. Furthermore, they have established a lower-bound for the regret for policy satisfying an appropriately defined consistency condition, and show that their policy was asymptotically efficient.
Later, the non-parametric context has been considered; several algorithms have been proposed, among which \emph{softmax action selection} policies and \emph{Upper-Confidence Bound} (UCB) policies.

Softmax methods are randomized policies where,  at time $t$, the arm $I_t$ is chosen at random by the player according to some probability distribution giving more weight to arms which have so-far performed well. The greedy action is given the highest selection probability, but all the others are ranked and weighted according to their accumulated rewards. The most common softmax action selection method uses a Gibbs, or Boltzman distribution. A prototypal example of softmax action selection is the so-called EXP3 policy (for \emph{Exponential-weight algorithm for Exploration and Exploitation}), which has been introduced by \cite{freund:schapire:robert:1995} for solving a worst-case sequential allocation problem and thouroughly examined as an instance of ``prediction with limited feedback'' problem in Chapter 6 of \cite{cesabianchi:lugosi:2006} (see also \cite{auer:cesabianchi:freund:schapire:2002,cesabianchi:lugosi:1999}).

UCB methods are deterministic policies extending the algorithm proposed by \cite{lai:robbins:1985} to a non-parametric context; they have been introduced and analyzed by \cite{agrawal:1995}. They consist in playing during the $t$-th round the arm $i$ that maximizes the upper bound of a confidence interval for expected reward $\mu(i)$, which is constructed from the past observed rewards.
The most popular, called UCB-1, relies on the upper-bound $\MReward{t}{i} +c_t(i)$, where $\MReward{t}{i}=(\Np[t][]{i})^{-1} \sum_{s=1}^t X_s(i)\1{I_s=i}$ denotes the empirical mean, and $c_t(i)$ is a \emph{padding function}.
A standard choice is $c_t(i)=B\sqrt{\xi\log (t)/\Np[t][]{i}}$, where $B$ is an upper-bound on the rewards and $\xi>0$ is some appropriate constant.
UCB-1 is defined in Algorithm \ref{algo:S-UCB}.
\begin{algorithm}
\caption{UCB-1}
\begin{algorithmic}
\label{algo:S-UCB}
\STATE for $t$ from $1$ to $K$, play arm $I_t=t$;
\STATE for $t$ from $K+1$ to $T$, play arm $$I_t=\argmax_{1\leq i\leq K} \MReward{t}{i} + c_t(i).$$
\end{algorithmic}
\end{algorithm}

UCB-1 belongs to the family of ``follow the perturbed leader'' algorithms, and has proven to retain the optimal logarithmic rate (but with suboptimal constant).
A finite-time analysis of this algorithm has been given in \cite{auer02finitetime,auer:2002, auer:cesabianchi:freund:schapire:2002}. Other types of padding functions are considered in \cite{audibert:munos:szepesvari:2007}.

\subsection{The non-stationary MAB problem}
In the non-stationary context, the rewards $\{X_s(i) \}_{s \geq 1}$ for arm $i$ are modeled by a sequence of independent random variables from potentially different distributions (unknown to the user) which may vary across time. For each $s > 0$, we denote by $\mu_s(i)$ the expectation of the reward $X_s(i)$ for arm $i$. Likewise, let $i^*_t$ be the arm with highest expected reward, denoted $\mu_t(*)$, at time $t$. The regret of a policy $\pi$ is now defined as the expected difference between the total rewards collected by the optimal policy $\pi^*$ (playing at each time instant  the arm $i^*_t$) and the total rewards collected by the policy $\pi$. Note that, in this paper, the non-stationary regret is not defined with respect to the best arm on average, but with respect to a strategy tracking the best arm at each step (this notion of regret is similar to the ``regret against arbitrary strategies'' introduced in Section 8 of \cite{auer:cesabianchi:freund:schapire:2002} for the non-stochastic bandit problem).

In this paper, we consider \emph{abruptly changing environments}: the distributions of rewards remain constant during periods and change at unknown time instants called \emph{breakpoints}.
In the following, we denote by $\Upsilon_T$  the number of abrupt changes in the reward distributions that occur before time $T$.
Another type of non-stationary MAB, where the distribution of rewards changes continuously, are considered in \cite{slivkins:upfal:2008}.

Standard soft-max and UCB policies are not appropriate for abruptly changing environments:
as stressed in \cite{hartland:gelly:baskiotis:teytaud:sebag:2006}, ``empirical evidence shows that
their Exploration versus Exploitation trade-off is not appropriate for abruptly changing
environments``.
To address this problem, several methods have been proposed.

In the family of softmax action selection policies, \cite{auer:cesabianchi:freund:schapire:2002}  and \cite{cesabianchi:lugosi:stoltz:2006, cesabianchi:lugosi:stoltz:2008} have proposed an adaptation referred to as \emph{EXP3.S} of the Fixed-Share algorithm, a computationally efficient variant of EXP3 called introduced by \cite{herbster:warmuth:1998} (see also \citep{cesabianchi:lugosi:2006} and the references therein).
Theorem 8.1 and Corollary 8.3  in \cite{auer:cesabianchi:freund:schapire:2002} state that when EXP3.S is tuned properly (which requires in particular that $\Upsilon_T$ is known in advance),
the expected regret is upper-bounded as
$$\E_\pi\left[R_T \right]\leq 2\sqrt{\rme-1}\sqrt{KT(\Upsilon_T\log(KT)+ \rme)}\eqsp.$$
Compared to the stationary case, such an upper-bound may seem deceiving: the rate $O(\sqrt{T\log T})$ is much larger than the $O(\log T)$ achievable in absence of changes. But actually, we prove in Section \ref{sec:lower-boundwithbreakpoint} that no policy can achieve an average regret smaller than $O(\sqrt{T})$ in the non-stationary case. Hence, EXP3.S matches the best achievable rate up to a factor $\sqrt{\log T}$. Moreover, by construction this algorithm can as well be used in an adversarial setup.

On the other hand, in the family of UCB policies, several attempts have been made; see for examples \cite{slivkins:upfal:2008} and \cite{kocsis:szepesvari:2006}.
In particular, \cite{kocsis:szepesvari:2006} have proposed an adaptation of the UCB policies that relies on a discount factor  $\gamma\in(0,1)$.
This policy constructs an UCB $\MReward{t}{\gamma,i}  + \Ic{i}{t}{\gamma}$
for the instantaneous expected reward, where the discounted empirical average is given by
$$\MReward{t}{\gamma,i}  =\frac{1}{\Nplay[][i]{t}{\gamma}}\sum_{s=1}^t \gamma^{t-s}X_s(i)\1{I_s=i}\eqsp, \quad \Nplay[][i]{t}{\gamma}=\sum_{s=1}^t \gamma^{t-s}\1{I_s=i},$$
and the discounted padding function is defined as
$$\Ic{i}{t}{\gamma}=2B\sqrt{\frac{\xi \log\Nplay[sum][]{t}{\gamma}}{\Nplay[][i]{t}{\gamma}}}\eqsp ,\quad  \Nplay[sum][]{t}{\gamma}=\sum_{i=1}^K \Nplay[][i]{t}{\gamma}\eqsp, $$
for an appropriate parameter $\xi$.
Using these notations, discounted-UCB (D-UCB) is defined in Algorithm \ref{algo:D-UCB}.
Remark that for $\gamma=1$, D-UCB boils down to the standard UCB-1 algorithm.

\begin{algorithm}
\caption{Discounted UCB}
\begin{algorithmic}\label{algo:D-UCB}
\STATE for $t$ from $1$ to $K$, play arm $I_t=t$;
\STATE for $t$ from $K+1$ to $T$, play arm $$I_t=\argmax_{1\leq i\leq K} \MReward{t}{\gamma,i} + \Ic{i}{t}{\gamma}.$$
\end{algorithmic}
\end{algorithm}

In order to estimate the instantaneous expected reward, the D-UCB policy averages past rewards with a discount factor giving more weight to recent observations.
We propose in this paper a more abrupt variant of UCB where averages are computed on a fixed-size horizon. At time $t$, instead of averaging the rewards over all past with a discount factor, \emph{sliding-window UCB} relies on a local empirical average of the observed rewards, using only the $\tau$ last plays.
Specifically, this algorithm constructs an UCB $\MReward[i]{t}{\tau}+ \Ic{i}{t}{\tau} $
for the instantaneous expected reward; the local empirical average is given by
$$\MReward[i]{t}{\tau}  = \frac{1}{\Nplay[][i]{t}{\tau}}\sum_{s=t-\tau+1}^tX_s(i)\1{I_s=i}\eqsp, \quad \Nplay[][i]{t}{\gamma}=\sum_{s=1}^t \gamma^{t-s}\1{I_s=i},$$
and the padding function is defined as
$$\Ic{i}{t}{\tau} = B\sqrt{\frac{\xi\log (t\wedge\tau)}{\Nplay[][i]{t}{\tau}}}\eqsp,$$
where $t\wedge\tau$ denotes the minimum of $t$ and $\tau$, and $\xi$ is a some appropriate constant.
The policy defined in Algorithm \ref{algo:SW-UCB} will be called in the sequel \emph{Sliding-Window UCB} (SW-UCB).

\begin{algorithm}
\caption{Sliding-Window UCB}
\begin{algorithmic}\label{algo:SW-UCB}
\STATE for $t$ from $1$ to $K$, play arm $I_t=t$;
\STATE for $t$ from $K+1$ to $T$, play arm $$I_t=\argmax_{1\leq i\leq K} \MReward[i]{t}{\tau} + \Ic{i}{t}{\tau},$$
\end{algorithmic}
\end{algorithm}

In this paper, we investigate the behaviors of the discounted-UCB and of the sliding-window-UCB in an abruptly changing environment, and prove that they are almost rate-optimal in a minimax sense.
In Section \ref{sec:discountedversionsofucb}, we derive a finite-time upper-bound on the regret of D-UCB.
In Section \ref{sec:slidingwindowucb}, we propose a similar analysis for the SW-UCB policy. We establish that it achieves the slightly better regret.
In Section \ref{sec:lower-boundwithbreakpoint}, we establish a lower-bound on the regret of any policy in an abruptly changing environment.
As a by-product, we show that any policy (like UCB-1) that achieves a logarithmic regret in the stationary case cannot reach a regret of order smaller than $T/\log T$ in presence of breakpoints.
The upper-bounds obtained in Sections \ref{sec:discountedversionsofucb} and \ref{sec:slidingwindowucb} are based on a novel deviation inequality for self-normalized averages with random number of summands which is stated and proved in Section \ref{sec:self-normalizedhoeffing}. A maximal inequality, of independent interest, is also derived in Section \ref{sec:self-normalizedmaxihoeffing}. Two simple Monte-Carlo experiments are presented to support our findings in Section \ref{sec:simulations}.

\section{Analysis of Discounted UCB}
\label{sec:discountedversionsofucb}
%
%
%
In this section, we analyze the behavior of D-UCB in an abruptly changing environment.
Let $\Upsilon_T$ denote the number of breakpoints before time $T$, and let $\Nbadplayed{i}{T}$ denote the number of times arm $i$ was played when it was not the best arm during the $T$ first rounds:
$$\Nbadplayed{i}{T}=\sum_{t=1}^T \1{I_t=i\neq i_t^*}.$$
Denote by $\Delta\mu_T(i)$ the minimum of the difference of expected reward of the best arm $\mu_t(*)$ and the expected reward $\mu_t(i)$ of the $i$-th arm   for all times $t\in\{1,\dots,T\}$ such that arm $i$ is not the leading arm ($i\neq i^*_t$),
\begin{equation}
\label{eq:definition-Delta}
\Delta\mu_T(i)= \min \left\{ t \in \{1, \dots, T\}, i \neq i^*_t, \mu_t(*) - \mu_t(i) \right\} \eqsp.
\end{equation}
We denote by $\P_{\DUCB{\gamma}}$ and $\E_{\DUCB{\gamma}}$ the probability distribution and expectation under policy D-UCB with discount factor $\gamma$. The next theorem computes a bound for the expected number of times in $T$ rounds that the arm $i$ is played, when this arm is suboptimal.

\begin{theorem}
Let $\xi>1/2$ and $\gamma\in(0,1)$.
For any arm $i\in\{1,\dots,K\}$,
\label{th:boundbadplayed}
\begin{equation}
\E_{\DUCB{\gamma}}\left[ \Nbadplayed{i}{T} \right]
\leq \operatorname{B}(\gamma) T (1-\gamma) \log \frac{1}{1-\gamma}\\ + \operatorname{C}(\gamma)\frac{\Upsilon_T}{1-\gamma}\log \frac{1}{1-\gamma}\eqsp,\label{eq:boundbadplayed}
\end{equation}
where
\begin{equation*}
\operatorname{B}(\gamma)=\frac{16B^2\xi }{\gamma^{1/(1-\gamma)}(\Delta\mu_T(i))^2}\frac{\lceil T(1-\gamma)\rceil}{T(1-\gamma)} \\+ \frac{2\left\lceil-\log(1-\gamma)/\log(1+4\sqrt{1-1/2\xi})\right\rceil }{-\log(1-\gamma)\left( 1-\gamma^{1/(1-\gamma)} \right)}
\end{equation*}
and
\begin{equation}
\label{eq:defCgamma}
\operatorname{C}(\gamma)=\frac{\gamma-1}{\log(1-\gamma)\log\gamma}\\\times \log\left((1-\gamma)\xi\log \Nplay[sum][]{K}{\gamma}\right)\eqsp.
\end{equation}
\end{theorem}
\begin{remark}
When $\gamma$ goes to $1$ we have $\operatorname{C}(\gamma) \to 1$ and
$$\operatorname{B}(\gamma)\to \frac{16 \rme B^2\xi}{(\Delta\mu_T(i))^2}+\frac{2}{(1-\rme^{-1})\log\left(1+4\sqrt{1- 1/2\xi}\right)}\eqsp.$$
\end{remark}

\begin{proof}
The proof is adapted from the finite-type analysis of \cite{auer02finitetime}.
There are however two main differences. First, because the expected reward changes, the discounted empirical mean
$\MReward[i]{t}{\gamma}$ is now a \emph{biased} estimator of the expected reward $\mu_t(i)$. The second difference stems from the deviation inequality itself: instead of using a Chernoff-Hoeffding bound, we use a novel tailored-made control on a self-normalized mean of the rewards with a random number of summands.
The proof is in 5 steps:
\paragraph{Step 1}
We upper-bound the number of times the suboptimal arm $i$ is played as follows:
\begin{align*}
\Nbadplayed{i}{T} &= 1+\sum_{t=K+1}^T \1{I_t=i\neq i_t^*}\nonumber\\
&= 1+\sum_{t=K+1}^T \1{I_t=i\neq i_t^*, \Nplay[][i]{t}{\gamma}<A(\gamma)}
+ \sum_{t=K+1}^T \1{I_t=i\neq i_t^*, \Nplay[][i]{t}{\gamma}\geq A(\gamma)} \eqsp,
\end{align*}
where
\begin{equation}
A(\gamma) = \frac{16B^2\xi\log \Nplay[sum][]{T}{\gamma}}{(\Delta\mu_T(i))^2}\eqsp.\label{eq:defA}
\end{equation}
Using Corollary \ref{cor:majabadondiscount} (stated and  proved in the Appendix), we may upper-bound the first sum in the RHS as:
\begin{equation*}
\sum_{t=K+1}^T \1{I_t=i\neq i_t^*, \Nplay[][i]{t}{\gamma}< A(\gamma)}
\leq \lceil T(1-\gamma)\rceil A(\gamma)\gamma^{-1/(1-\gamma)} \eqsp. 
\end{equation*}
In the sequel, for any positive $m$, we denote by
$\T(\gamma)$ the set of all indices $t\in \{K+1,\dots, T\}$ such that $\mu_s(j)=\mu_t(j)$ for all $j\in\{1,\dots,K\}$ and all $t-D(\gamma)<s\leq t$, where $$D(\gamma)=\frac{\log\left((1-\gamma)\xi\log \Nplay[sum][]{K}{\gamma} \right)}{\log \gamma}\eqsp.$$
During a number of rounds (that depends on $\gamma$) following a breakpoint, the estimates of the expected rewards can be  poor. 
Because of this, the D-UCB policy may play constantly the suboptimal arm $i$, which leads to the following bound:
\begin{equation*}
\sum_{t=K+1}^T \1{I_t=i\neq i_t^*,  \Nplay[][i]{t}{\gamma}\geq A(\gamma)}
\leq  \Upsilon_TD(\gamma)
+ \sum_{t\in\T(\gamma)} \1{I_t=i\neq i_t^*,  \Nplay[][i]{t}{\gamma}\geq A(\gamma)}\eqsp.
\end{equation*}
Putting everything together, we obtain:
\begin{equation}
\Nbadplayed{i}{T} \leq 1+ \lceil T(1-\gamma)\rceil A(\gamma)\gamma^{-1/(1-\gamma)} +\Upsilon_TD(\gamma)
+ \sum_{t\in\T(\gamma)} \1{I_t=i\neq i_t^*,  \Nplay[][i]{t}{\gamma}\geq A(\gamma)}\eqsp.
\label{eq:ducb:pet}
\end{equation}

\paragraph{Step 2}
Now, for $t\in\T(\gamma)$  the event $\left\{I_t=i\neq i^*_t, \Nplay[][i]{t}{\gamma}\geq A(\gamma)\right\}$ may be decomposed as follows:
\begin{multline}
\left\{I_t=i\neq i^*_t, \Nplay[][i]{t}{\gamma}\geq A(\gamma)\right\}
\subseteq \left\{\MReward[i]{t}{\gamma}>\mu_t(i) + \Ic{i}{t}{\gamma}\right\}
\cup \left\{\MReward[*]{t}{\gamma}<\mu_t(*) - \Ic{*}{t}{\gamma}\right\} \\
\cup \left\{\mu_t(*)-\mu_t(i) < 2\Ic{i}{t}{\gamma}, \Nplay[][i]{t}{\gamma}\geq A(\gamma)\right\}. \label{eq:ducb:decomp}
\end{multline}
In words, playing the suboptimal arm $i$ at time $t$ may occur in three cases: if $\mu_t(i)$ is substantially over-estimated, if $\mu_t(*)$ is substantially under-estimated, or if $\mu_t(i)$ and $\mu_t(*)$ are close from each other.
But for the choice of $A(\gamma)$ given in Equation (\ref{eq:defA}), we have
\begin{align*}
\Ic{i}{t}{\gamma} 
&\leq  2B\sqrt{\frac{\xi\log \Nplay[sum][]{t}{\gamma}}{A(\gamma)}} \leq \frac{\Delta\mu_T(i)}{2}\eqsp,
\end{align*}
so that the event $\left\{\mu_t(*)-\mu_t(i) < 2\Ic{i}{t}{\gamma},  \Nplay[][i]{t}{\gamma}\geq A(\gamma)\right\}$ never occurs.

In Steps 3 and 4 we upper-bound the probability of the two first events of the RHS of \eqref{eq:ducb:decomp}.
We show that for $t\in\T(\gamma)$, that is at least $D(\gamma)$ rounds after a breakpoint, the expected rewards of all arms are well estimated with high probability. For all $j\in\{1, \dots, K\}$, consider the following events
\begin{equation*}
\mathcal{E}_t(\gamma, j)=\big\{\MReward[i]{t}{\gamma}>\mu_t(j) + \Ic{j}{t}{\gamma}\big\}
\end{equation*}
The idea is the following: we upper-bound the probability of $\mathcal{E}_t(\gamma, j)$ 
by separately considering the fluctuations of $\MReward[j]{t}{\gamma}$ around $\EReward{t}{\gamma, j} /\Nplay[][j]{t}{\gamma}$, and the `bias' $\EReward{t}{\gamma, j}/\Nplay[][j]{t}{\gamma} -\mu_t(j)$, where
$$\EReward{t}{\gamma, j}=\sum_{s=1}^{t} \gamma^{t-s}\1{I_s=j}\mu_s(j)\eqsp.$$

\paragraph{Step 3}
Let us first consider the bias.
First note that $\EReward{t}{\gamma, j}/\Nplay[][j]{t}{\gamma}$, as a convex combination of elements $\mu_s(j)\in[0,B]$, belongs to interval $[0,B]$.
Hence, $|\EReward{t}{\gamma, j}/\Nplay[][j]{t}{\gamma}-\mu_t(j)|\leq B$.
Second, for $t\in\T(\gamma)$,
\begin{multline*}
\left| \EReward{t}{\gamma, j} - \mu_t(j) \Nplay{t}{\gamma}\right|
= \left|\sum_{s=1}^{t-D(\gamma)} \gamma^{t-s}\left( \mu_s(j)-\mu_t(j)\right)\1{I_s=j}  \right|\\
\leq \sum_{s=1}^{t-D(\gamma)} \gamma^{t-s}\left|\mu_s(j)-\mu_t(j)\right|\1{I_s=j} \leq B\sum_{s=1}^{t-D(\gamma)} \gamma^{t-s}\1{I_s=j}
 = B\gamma^{D(\gamma)} \Nplay[][j]{t-D(\gamma)}{\gamma}.
\end{multline*}
As $ \Nplay[][j]{t-D(\gamma)}{\gamma}\leq (1-\gamma)^{-1}$, we get
$|\EReward{t}{\gamma, j}/\Nplay[][j]{t}{\gamma}-\mu_t(j)|\leq B \gamma^{D(\gamma)}(1-\gamma)^{-1}$.
Altogether,
$$\left|\frac{\EReward{t}{\gamma, j}}{\Nplay[][j]{t}{\gamma}}-\mu_t(j)\right|\leq B \left(1\wedge \gamma^{D(\gamma)}(1-\gamma)^{-1}  \right).$$
Hence, using the elementary inequality $1\wedge x\leq \sqrt{x}$ and the definition of $D(\gamma)$, we obtain for $t\in\T(\gamma)$:
$$\left|\frac{\EReward{t}{\gamma, j}}{\Nplay[][j]{t}{\gamma}}-\mu_t(j)\right|\leq B\sqrt{ \frac{\gamma^{D(\gamma)}}{(1-\gamma)\Nplay[][i]{t}{\gamma}}} \leq B\sqrt{\frac{\xi\log\Nplay[sum][]{K}{\gamma}}{\Nplay[][j]{t}{\gamma}}} \leq \frac{1}{2}\Ic{j}{t}{\gamma}\eqsp.$$
In words:  $D(\gamma)$ rounds after a breakpoint, the `bias' is smaller than the half of the padding function.
The other half of the padding function is used to control the fluctuations.
In fact, for $t\in\T(\gamma)$: 
\begin{align*}
\P_{\DUCB{\gamma}}\left( \mathcal{E}_t(\gamma, j)
 \right)
&\leq \P_{\DUCB{\gamma}}\Bigg( \MReward[j]{t}{\gamma}>\mu_t(j)+ B\sqrt{\frac{\xi\log \Nplay[sum][]{t}{\gamma}}{\Nplay[][j]{t}{\gamma}}}+\left|\frac{\EReward{t}{\gamma, j}}{\Nplay[][j]{t}{\gamma}}-\mu_t(j)\right|\Bigg)\\
&\leq \P_{\DUCB{\gamma}}\Bigg( \MReward[j]{t}{\gamma}-\frac{\EReward{t}{\gamma, j}}{\Nplay[][j]{t}{\gamma}}> B\sqrt{\frac{\xi\log \Nplay[sum][]{t}{\gamma}}{\Nplay[][j]{t}{\gamma}}}\Bigg)\eqsp.
\end{align*}
\paragraph{Step 4}
Denote the discounted total reward obtained with arm $j$ by $$\Reward{t}{\gamma, j}=\sum_{s=1}^t \gamma^{t-s}\1{I_s=j}X_s(j)=\Nplay[][j]{t}{\gamma}\MReward[j]{t}{\gamma}\eqsp.$$
Using Theorem \ref{th:concentrationR} and the fact that $\Nplay[][j]{t}{\gamma}\geq \Nplay[][j]{t}{\gamma^2}$, the previous inequality rewrites:
\begin{align*}
\P_{\DUCB{\gamma}}\left( \mathcal{E}_t(\gamma, j) \right) &\leq
 \P_{\DUCB{\gamma}}\Bigg( \frac{\Reward{t}{\gamma, j}-\EReward{t}{\gamma, j}}{\sqrt{\Nplay[][j]{t}{\gamma^2}}}> B\sqrt{\frac{\xi \Nplay[][j]{t}{\gamma}\log \Nplay[sum][]{t}{\gamma}}{\Nplay[][j]{t}{\gamma^2}}}\Bigg)\nonumber\\
&\leq \P_{\DUCB{\gamma}}\Bigg( \frac{\Reward{t}{\gamma, j}-\EReward{t}{\gamma, j}}{\sqrt{\Nplay[][j]{t}{\gamma^2}}}> B\sqrt{\xi\log \Nplay[sum][]{t}{\gamma}}\Bigg)\nonumber\\
&\leq
\left\lceil\frac{\log \Nplay[sum][]{t}{\gamma}}{\log(1+\eta)}\right\rceil \exp\left(-2\xi\log\Nplay[sum][]{t}{\gamma} \left(1-\frac{\eta^2}{16}\right)\right)\nonumber\\
&=\left\lceil\frac{\log \Nplay[sum][]{t}{\gamma}}{\log(1+\eta)}\right\rceil \Nplay[sum][]{t}{\gamma}^{-2\xi\left(1-\frac{\eta^2}{16}\right)}\eqsp.
\end{align*}



\paragraph{Step 5}
Hence, we finally obtain from Equation \eqref{eq:ducb:pet} :
\begin{equation*}
\E_{\DUCB{\gamma}}\left[ \Nbadplayed{i}{T} \right]\leq 1+\lceil T(1-\gamma)\rceil A(\gamma)\gamma^{-1/(1-\gamma)} +  D(\gamma) \Upsilon_T \\
 +2\sum_{t\in\T(\gamma)} \left\lceil\frac{\log \Nplay[sum][]{t}{\gamma}}{\log(1+\eta)}\right\rceil \Nplay[sum][]{t}{\gamma}^{-2\xi\left(1-\frac{\eta^2}{16}\right)}\eqsp.
\end{equation*}
When $\Upsilon_T\neq 0$, $\gamma$ is taken strictly smaller than $1$ (see Remark \ref{rk:gammaofRn}).
As $\xi>\frac{1}{2}$, we take $\eta=4\sqrt{1- 1/2\xi}$, so that $2\xi\left(1-\eta^2/16\right) =1$.
For that choice, with  $\tau=(1-\gamma)^{-1}$,
\begin{align*}
\sum_{t\in\T(\gamma)}\left\lceil\frac{\log \Nplay[sum][]{t}{\gamma}}{\log(1+\eta)}\right\rceil \Nplay[sum][]{t}{\gamma}^{-2\xi\left(1-\frac{\eta^2}{16}\right)} &\leq \tau -K + \sum_{t=\tau}^T\left\lceil\frac{\log \Nplay[sum][]{\tau}{\gamma}}{\log(1+\eta)}\right\rceil \Nplay[sum][]{\tau}{\gamma}^{-1}\\
& \leq \tau -K + \left\lceil\frac{\log \Nplay[sum][]{\tau}{\gamma}}{\log(1+\eta)}\right\rceil\frac{n}{\Nplay[sum][]{\tau}{\gamma}}\\
&\leq \tau -K +  \left\lceil\frac{\log\frac{1}{1-\gamma}}{\log(1+\eta)}\right\rceil\frac{T(1-\gamma)}{1-\gamma^{1/(1-\gamma)}},
\end{align*}
we obtain the statement of the Theorem.
\end{proof}

\begin{remark}
\label{rk:gammaofRn}
If horizon $T$ and the growth rate of the number of breakpoints $\Upsilon_T$ are known in advance, the discount factor $\gamma$ can be chosen so as to minimize the RHS in Equation \ref{eq:boundbadplayed}.
Taking $\gamma=1-(4B)^{-1}\sqrt{\Upsilon_T/T}$ yields:
$$\E_{\DUCB{\gamma}}\left[ \Nbadplayed{i}{T} \right]= O\left( \sqrt{T \Upsilon_T}\log T \right)\eqsp .$$
Assuming that $\Upsilon_T=O(T^\beta)$ for some $\beta\in[0,1)$, the regret is upper-bounded as
$O\left( T^{(1+\beta)/2}\log T \right)$.
In particular, if $\beta=0$, the number of breakpoints $\Upsilon_T$ is upper-bounded by $\Upsilon$ independently of $T$, taking $\gamma=1-(4B)^{-1}\sqrt{\Upsilon/T}$, the regret is bounded by $O\left( \sqrt{\Upsilon T}\log T \right)$.
Thus, D-UCB matches the lower-bound of Theorem \ref{th:binf} up to a factor $\log T$.
\end{remark}

\begin{remark}
On the other hand, if the breakpoints have a positive density over time (say, if $\Upsilon_T\leq rT$ for a small positive constant $r$), then $\gamma$ has to remain lower-bounded independently of $T$; Theorem \ref{th:boundbadplayed} gives a linear, non-trivial bound on the regret and permits to calibrate the discount factor $\gamma$ as a function of the density of the breakpoint: taking $\gamma=1-\sqrt{r}/(4B)$ we get an upper-bound with a dominant term in 
$O\left(-T\sqrt{r}\log r\right)$.
\end{remark}


\begin{remark}
Theorem \ref{th:concentrationRmaxi} shows that for $\xi>1/2$ and $t\in\T(\gamma)$, with high probability $\MReward[i]{t}{\gamma}$ is actually \emph{never} larger than $\mu_t(i)+\Ic{i}{t}{\gamma}$.
\end{remark}

\begin{remark}
If the growth rate of $\Upsilon_T$ is known in advance, but not the horizon $T$, then we can use the ``doubling trick'' to set the value of $\gamma$.
Namely, for $t$ and $k$ such that $2^k\leq t<2^{k+1}$, take $\gamma=1-(4B)^{-1}(2^k)^{(\beta-1)/2}$.
\end{remark}


\section{Sliding window UCB}
\label{sec:slidingwindowucb}
In this section, we analyze the performance of SW-UCB in an abruptly changing environment.
We denote by $\P_{\SWUCB{\tau}}$ and $\E_{\SWUCB{\tau}}$ the probability distribution and expectation under policy SW-UCB with window size $\tau$.
\begin{theorem}
\label{th:boundbadplayedwindow}
Let $\xi>1/2$.
For any integer $\tau$ and any arm $i\in\{1,\dots,K\}$,
\begin{equation}
\E_{\SWUCB{\tau}}\left[ \Nbadplayed{i}{T} \right]\leq
\operatorname{C}(\tau)\frac{T\log\tau}{\tau} + \tau \Upsilon_T +\log^2(\tau)\eqsp,
\label{eq:boundbadplayedwindow}
\end{equation}
where
\begin{equation*}
\operatorname{C}(\tau)=\frac{4B^2\xi}{(\Delta\mu_T(i))^2}\frac{\lceil T/\tau\rceil}{T/\tau} \\+ \frac{2}{\log\tau}\left\lceil\frac{\log(\tau)} {\log(1+4\sqrt{1-(2\xi)^{-1}})}\right\rceil\eqsp.
\end{equation*}
\end{theorem}
\begin{remark}
As $\tau$ goes to infinity
$$\operatorname{C}(\tau)\to \frac{4B^2\xi}{(\Delta\mu_T(i))^2}+\frac{2}{\log(1+4\sqrt{1-(2\xi)^{-1}})}.$$
\end{remark}

\begin{proof}
We follow the lines of the proof of Theorem \ref{th:boundbadplayed}. The main difference is that for $t\in\T(\tau)$ defined here as the set of all indices $t\in \{K+1,\dots, T\}$ such that $\mu_s(j)=\mu_t(j)$ for all $j\in\{1,\dots,K\}$ and all $t-\tau<s\leq t$, the bias exactly vanishes; consequently, Step 3 can be bypassed.
\paragraph{Step 1}
Let $A(\tau) = 4B^2\xi\log \tau(\Delta\mu_T(i))^{-2}$; using Lemma \ref{lem:majabadondiscount}, we have:
\begin{align}
\Nbadplayed{i}{T}& = 1+\sum_{t=K+1}^T \1{I_t=i\neq i_t^*} \nonumber\\
&\leq 1+\sum_{t=1}^T \1{I_t=i, \Nplay[][i]{t}{\tau}< A(\tau)}+ \sum_{t=K+1}^T \1{I_t=i\neq i_t^*, \Nplay[][i]{t}{\tau}\geq A(\tau)}\nonumber\\
&\leq 1+\lceil{T/\tau}\rceil A(\tau)+ \sum_{t=K+1}^T \1{I_t=i\neq i_t^*, \Nplay[][i]{t}{\tau}\geq A(\tau)}\nonumber\\
&\leq 1+\lceil{T/\tau}\rceil A(\tau) +  \Upsilon_T\tau+ \sum_{t\in\T(\tau)} \1{I_t=i\neq i_t^*, \Nplay[][i]{t}{\tau}\geq A(\tau)}\label{eq:majnbadplayedsw}
\end{align}
\paragraph{Step 2}
For $t\in\T(\tau)$ we have
\begin{multline}
\left\{I_t=i, \Nplay[][i]{t}{\tau}\geq A(\tau)\right\}\subset\left\{\MReward[i]{t}{\tau}>\mu_t(i) + \Ic{i}{t}{\tau}\right\}
\cup \left\{\MReward[*]{t}{\tau}<\mu_t(*) - \Ic{*}{t}{\tau}\right\}\\
\cup \left\{\mu_t(*)-\mu_t(i) < 2\Ic{i}{t}{\tau}, \Nplay[][i]{t}{\tau}\geq A(\tau)\right\}. \label{eq:decompsw}
\end{multline}
On the event $\left\{\Nplay{t}{\tau, i}\geq A(\tau)\right\}$, we have
\begin{equation*}
\Ic{i}{t}{\tau} =B\sqrt{\frac{\xi\log(t\wedge\tau)}{\Nplay[][i]{t}{\tau}}}\\
\leq B\sqrt{\frac{\xi\log \tau}{A(\tau)}}
= B\sqrt{\frac{\xi\log(\tau) \, (\Delta\mu_T(i))^2}{4B^2\xi\log \tau}}
\leq \frac{\Delta\mu_T(i)}{2},
\end{equation*}
so that the event $\left\{\mu_t(*)-\mu_t(i) < 2\Ic{i}{t}{\tau}, \Nplay[][i]{t}{\tau}\geq A(\tau)\right\}$ has $\P_{\SWUCB{\tau}}$-probability $0$.
\paragraph{Steps 3-4}
Now, for $t\in\T(\tau)$ and for all $j\in\{1,\dots,K\}$, Corollary \ref{cor:concentrationMR} applies and yields:
\begin{align}
\P\left( \MReward[j]{t}{\tau}>\mu_t(j) + \Ic{j}{t}{\tau}\right)
&\leq \P\bigg( \MReward[j]{t}{\tau}>\mu_t(j) + B\sqrt{\frac{\xi\log (t\wedge\tau)}{\Nplay[][j]{t}{\tau}}}\Bigg)\nonumber\\
&\leq \left\lceil\frac{\log(t\wedge\tau)}{\log(1+\eta)}\right\rceil \exp\left(-2\xi\log (t\wedge\tau)\left(1-\frac{\eta^2}{16}\right)\right)\nonumber\\
&= \left\lceil\frac{\log (t\wedge\tau)}{\log(1+\eta)}\right\rceil (t\wedge\tau)^{-2\xi\left(1-\eta^2/16\right)}\eqsp,\label{eq:majproba1sw}
\end{align}
and similarly
\begin{equation}
\P\left( \MReward[j]{t}{\tau}<\mu_t(j) - \Ic{j}{t}{\tau}\right) \\\leq \left\lceil\frac{\log(t\wedge\tau)}{\log(1+\eta)}\right\rceil (t\wedge\tau)^{-2\xi\left(1-\eta^2/16\right)}.
\label{eq:majproba2sw}
\end{equation}
\paragraph{Steps 5}
In the following we take $\eta=4\sqrt{1-\frac{1}{2\xi}}$, so that we have $2\xi\left(1-\eta^2/16 \right)=1$.
Thus, using Equations \eqref{eq:decompsw},\eqref{eq:majproba1sw} and \eqref{eq:majproba2sw}, Inequality \eqref{eq:majnbadplayedsw} yields
\begin{align*}
\E_{\SWUCB{\tau}}\left[ \Nbadplayed{i}{T} \right]&\leq 1+\lceil T/\tau \rceil A(\tau) + \tau \Upsilon_T + 2\sum_{t=1}^T \frac{\left\lceil\frac{\log(t\wedge\tau)} {\log(1+\eta)}\right\rceil}{(t\wedge\tau)}.
\end{align*}
The results follows, noting that
\begin{equation*}
\sum_{t=K+1}^T \frac{\log(t\wedge\tau)}{t\wedge\tau} \leq \sum_{t=2}^\tau \frac{\log t}{t} + \sum_{t=1}^T \frac{\log \tau}{\tau}\\
 \leq \frac{1}{2}\log^2(\tau) + \frac{T\log \tau}{\tau}.
\end{equation*}

\end{proof}

\begin{remark}
 If the horizon $T$ and the growth rate of the number of breakpoints $\Upsilon_T$ are known in advance, the window size $\tau$ can be chosen so as to minimize the RHS in Equation \eqref{eq:boundbadplayedwindow}.
Taking $\tau=2B\sqrt{T\log(T)/\Upsilon_T}$ yields
$$\E_{\SWUCB{\tau}}\left[ \Nbadplayed{i}{T} \right]= O\left( \sqrt{\Upsilon_T T\log T} \right)\eqsp.$$
Assuming that $\Upsilon_T=O(T^\beta)$ for some $\beta\in[0,1)$, the average regret is upper-bounded as
$O\left( T^{(1+\beta)/2}\sqrt{\log T} \right).$
In particular, if $\beta=0$, the number of breakpoints $\Upsilon_T$ is upper-bounded by $\Upsilon$ independently of $T$, then with $\tau=2B\sqrt{T\log(T)/\Upsilon }$ the upper-bound is $O\left( \sqrt{\Upsilon T\log T} \right)$.
Thus, SW-UCB matches the lower-bound of Theorem \ref{th:binf} up to a factor $\sqrt{\log T}$, slightly better than the D-UCB.
\end{remark}
\begin{remark}
On the other hand, if the breakpoints have a positive density over time, then $\tau$ has to remain lower-bounded independently of $T$.
For instance, if $\Upsilon_T\leq rT$ for some (small) positive rate $r$, and for the choice
$ \tau=2B\sqrt{-\log r /r}$,
Theorem \ref{th:boundbadplayedwindow} gives
$$\E_{\SWUCB{\tau}}\left[ \Nbadplayed{i}{T} \right] = O\left( T\sqrt{-r\log\left(r\right)}  \right).$$
\end{remark}

\begin{remark}
\label{rem:req0window}
If there is no breakpoint ($\Upsilon_T=0$), the best choice is obviously to take the window as a large as possible, that is $\tau=T$.
Then the procedure is exactly standard UCB.
A slight modification of the preceeding proof for $\xi=\frac{1}{2}+\epsilon$ with arbitrary small $\epsilon$ yields
\begin{align*}
\E_{\mathrm{UCB}}\left[ \Nbadplayed{i}{T} \right]&\leq \frac{2B^2}{(\Delta\mu(i))^2}\log T\left( 1+O(1) \right).
\end{align*}
We recover the same kind of bounds that are usually obtained in the analysis of UCB, see for instance \cite{auer02finitetime}, with a better constant.
\end{remark}

\begin{remark}
 The computational complexity of SW-UCB is, as for D-UCB, linear in time and does not involve $\tau$. However, SW-UCB requires to store the last $\tau$ actions and rewards at each time $t$ in order to efficiently update $\Nplay[][i]{t}{\tau}$ and $\MReward[i]{t}{\tau}$.
\end{remark}

\section{A lower-bound on the regret in abruptly changing environment}
\label{sec:lower-boundwithbreakpoint}
In this section, we consider a particular non-stationary bandit problem where the distributions of rewards on each arm are piecewise constant and have two breakpoints. Given any policy $\pi$, we derive a lower-bound on the number of times a sub-optimal arm is played (and thus, on the regret) in at least one such game.
Quite intuitively, the less explorative a policy is, the longer it may keep a suboptimal policy after a breakpoint.
Theorem \ref{th:binf} gives a precise content to this statement.


As in the previous section, $K$ denotes the number of arms, and the rewards are assumed to be bounded in $[0,B]$.
Consider any deterministic policy $\pi$ of choosing the arms  $I_1, \ldots, I_T$ played at each time depending to the past rewards $$G_t\eqdef X_{t}(I_t),$$ and recall that $I_t$ is measurable with respect to the sigma-field $\sigma(G_1, \dots, G_t)$ of the past observed rewards.
Denote by $\Np[s:t]{i}$ the number of times arm $i$ is played between times $s$ and $t$
$$\Np[s:t]{i} = \sum_{u=s}^t \1{I_u=i},$$
and $\Np[T]{i}=\Np[1:T]{i}$.
For $1\leq i\leq K$, let $\ps_i$ be the probability distribution of the outcomes of arm $i$, and let $\mus(i)$ denote its expectation.
Assume that $\mus(1)>\mus(i)$ for all $2\leq i\leq K$.
Denote by $\Ps$ the distribution of rewards under policy $\pi$, that is:
$$d\Ps(g_{1:T} | I_{1:T})=\prod_{t=1}^T d\ps_{i_t}(g_t).$$
For any random variable $W$ measurable with respect to $\sigma(G_1,\dots,G_T)$, denote by $\Es[W]$ its expectation under distribution $\Ps$.

In the sequel, we divide the period $\{1,\ldots,T\}$ into epochs of size $\tau\in \{1,\dots, T\}$, and we modify the distribution of the rewards so that on one of those periods, arm $K$ becomes the one with highest expected reward.
Specifically: let $Q$  be a distribution of rewards with expectation $\nu>\mus(1)$, let $\delta=\nu-\mus(1)$ and let $\alpha=\KL(\ps_K; Q)$ be the Kullback-Leibler divergence between $\ps_K$ and $Q$.
For all $1\leq j\leq M=\left\lfloor \frac{T}{\tau}\right\rfloor$, we consider the modification $\P_\pi^j$ of $\Ps$ such that on the $j$-th period of size $\tau$, the distribution of rewards of the $K$-th arm is changed to $\nu$.
That is, for every sequence of rewards $g_{1:T}$,
$$\frac{d\P_\pi^j}{d\Ps}(g_{1:T} | I_{1:T})  =\prod_{t=1+(j-1)\tau, I_t=K}^{j\tau} \frac{d Q}{d\ps_K} (g_t) \eqsp.$$
Besides, let
$$\Np[][j]{i}=\Np[1+(j-1)\tau : j\tau][]{i}$$
be the number of times arm $i$ is played in the $j$-th period.
For any random variable $W$ in $\sigma(G_1,\dots,G_T)$, denote by $ \E_\pi^j[W]$ its expectation under distribution $\P_\pi^j$.
Now, denote by $\P_\pi^*$ the distribution of rewards when $j$ is chosen uniformly at random in the set $\{1,\dots, M\}$ - in other words, $\P_\pi^*$ is the (uniform) mixture of the $(\P_\pi^j)_{1\leq j\leq M}$, and denote by $\E_\pi^*[\cdot]$ the expectation under $\P_\pi^*$:
$$\E_\pi^*[W]=\frac{1}{M}\sum_{j=1}^M\E_\pi^j[W].$$
In the following, we lower-bound the expect regret of any policy $\pi$ under $\P_\pi^*$ in terms of its regret under $\P_\pi$.
\begin{theorem}
\label{th:binf}
For any policy $\pi$ and any horizon $T$ such that $64/(9\alpha)\leq \Es[\Np[T][]{K}]\leq T/(4\alpha)$,
\begin{equation*}
 \E_\pi^*[R_T] \geq C(\mu)\frac{T}{\E_\pi[R_T]},
\end{equation*}
where 
$$C(\mu)=\frac{32\delta(\mu(1)-\mu(K))}{27\alpha}\eqsp.$$
\end{theorem}
\begin{proof}
The main ingredients of this reasoning are inspired by the proof of Theorem 5.1 in \cite{auer:cesabianchi:freund:schapire:2002}, see also \cite{kulkarni:lugosi:2000}.
First, note that the K\"ullback-Leibler divergence $\KL(\Ps, \P_\pi^j)$ is:
\begin{eqnarray*}
\KL(\Ps, \P_\pi^j) & =&\sum_{t=1}^T \KL\left( \Ps\left( G_t | G_{1:t-1}\right) ; \P_\pi^j\left( G_t | G_{1:t-1}\right)\right)\\
& = & \sum_{t=1+(j-1)\tau}^{j\tau} \Ps\left( I_t = K\right) \KL(\ps_K; Q) \\
& = & \alpha \Es\left[\Np[1+(j-1)\tau: j\tau]{K}\right] \eqsp.
\end{eqnarray*}
Hence, by Lemma A.1 in \cite{auer:cesabianchi:freund:schapire:2002},
\begin{eqnarray*}
 \E_\pi^j[\Np[][j]{K}] \leq \Es[\Np[][j]{K}] + \frac{\tau}{2}\sqrt{ \KL(\Ps, \P_\pi^j)}
=\Es[\Np[][j]{K}] + \frac{\tau}{2}\sqrt{ \alpha \Es[\Np[][j]{K}]} \eqsp.
\end{eqnarray*}
Consequently, since $\sum_{j=1}^M  \Np[][j]{K} \leq \Np[T][]{K}$,
\begin{align*}
\sum_{j=1}^M  \E_\pi^j[\Np[][j]{K}]
\leq \E[\Np[T][]{K}] + \frac{\tau}{2} \sum_{j=1}^M \sqrt{\alpha \Es[\Np[][j]{K}]}
\leq  \Es[\Np[T][]{K}] + \frac{\tau}{2} \sqrt{\alpha M\Es[\Np[T][]{K}]} \eqsp.
\end{align*}
Thus, there exists $1\leq j\leq M$ such that
\begin{align*}
  \E_\pi^*[\Np[][j]{K}] 
&\leq \frac{1}{M}\Es[\Np[T][]{K}] + \frac{\tau}{2M} \sqrt{\alpha M\Es[\Np[T][]{K}]}\\
& \leq  \frac{\tau}{T-\tau}\Es[\Np[T][]{K}]+ \frac{1}{2} \sqrt{\alpha \frac{\tau^3}{T-\tau}\Es[\Np[T][]{K}]} \eqsp.
\end{align*}
Now, the expectation under $\P_\pi^*$ of the regret $R_T$ is lower-bounded by:
\begin{equation*}
\E_\pi^*[R_T] 
\geq  \delta\left(\tau -  \E_\pi^*[\Np[T][]{K}]\right)\\
\geq  \delta \left(\tau - \frac{\tau}{T-\tau}\Es[\Np[T][]{K}] - \frac{1}{2} \sqrt{\alpha\frac{\tau^3}{T-\tau}\Es[\Np[T][]{K}]}\right) \eqsp.
\end{equation*}
Maximizing the right hand side of the previous inequality by choosing $\tau= 16T/(9\alpha\Es[\Np{K}])$ or equivalently $M=9\alpha/(16\E_\pi[N(K)])$ leads to the lower-bound:
\begin{equation*}
 \E_\pi^*[R_T] \geq \frac{32\delta}{27\alpha}\left( 1-\frac{\alpha\Es[\Np[T][]{K}]}{T} \right)^2\left( 1-\frac{16}{9\alpha\Es[\Np[T][]{K}]} \right)\frac{T}{\Es[\Np[T][]{K}]} \eqsp.
\end{equation*}
To conclude, simply note that $\Np[T][]{K}\leq \E_\pi[R_T]/(\mu(1)-\mu(K))$. We obtain:
\begin{equation*}
 \E_\pi^*[R_T] \geq \frac{32\delta(\mu(1)-\mu(K))}{27\alpha}\left( 1-\frac{\alpha\Es[\Np[T][]{K}]}{T} \right)^2\\ \left( 1-\frac{16}{9\alpha\Es[\Np[T][]{K}]} \right)\frac{T}{\E_\pi[R_T]} \eqsp,
\end{equation*}
which directly leads to the statement of the Theorem.
\end{proof}

The following corollary states that no policy can have a non-stationary regret of order smaller than $\sqrt{T}$.
It appears here as a consequence of Theorem \ref{th:binf}, although it can also be proved directly.
\begin{corollary}
For any policy $\pi$ and any positive horizon $T$,
$$\max\{\E_\pi(R_T), \E_\pi^*(R_T)\} \geq \sqrt{C(\mu)T}\eqsp.$$
\end{corollary}
\begin{proof}
If $\Es[\Np[T][]{K}]\leq 16/(9\alpha)$, or if $\Es[\Np[T][]{K}] \geq T/\alpha$, the result is obvious. Otherwise, Theorem \ref{th:binf} implies that:
$$\max\{\E_\pi(R_T), \E_\pi^*(R_T)\} \geq \max\{\E_\pi(R_T), C(\mu)\frac{T}{\E_\pi(R_T)}\} \geq \sqrt{C(\mu)T} \eqsp.$$
\end{proof}

\begin{remark}
To keep simple notations, Theorem \ref{th:binf} is stated and proved here for deterministic policy.
It is easily verified that the same results also holds for randomized strategies (such as EXP3-P, see \cite{auer:cesabianchi:freund:schapire:2002}).
\end{remark}
\begin{remark}
In words, Theorem \ref{th:binf} states that for any policy not playing each arm often enough, there is necessarily a time where a breakpoint is not seen after a long period.
For instance, as standard UCB satisfies $\Es[\Np{K}] = \Theta(\log{T})$, then
\begin{eqnarray*}
 \E_\pi^*[R_T] & \geq & c \frac{T}{\log T}
\end{eqnarray*}
for some positive $c$ depending on the reward distribution.
\end{remark}

\begin{remark}
This result is to be compared with standard minimax lower-bounds on the regret.
On one hand, a \emph{fixed-game lower-bound} in $O(\log T)$ was proved in \cite{lai:robbins:1985} for the stationary case, when the distributions of rewards are fixed and $T$ is allowed to go to infinity.
On the other hand, a finite-time \emph{minimax lower-bound} for individual sequences  in $O(\sqrt{T})$ is proved in \cite{auer:cesabianchi:freund:schapire:2002}. In this bound, for each horizon $T$ the worst case among all possible reward distributions is considered, which explains the discrepancy.
This result is obtained by letting the distance between distributions of rewards tend to 0 (typically, as $1/\sqrt{T}$).
In Theorem \ref{th:binf}, no assumption is made on the distributions of rewards $\ps_i$ and $Q$, their distance actually remains lower-bounded independently of $T$.
In fact, in the case considered here minimax regret and fixed-game minimal regret appear to have the same order of magnitude.
\end{remark}

\section{Simulations}
\label{sec:simulations}
\begin{figure}
\centering
\begin{tabular}{cc}
\includegraphics[width=7cm,height=2cm]{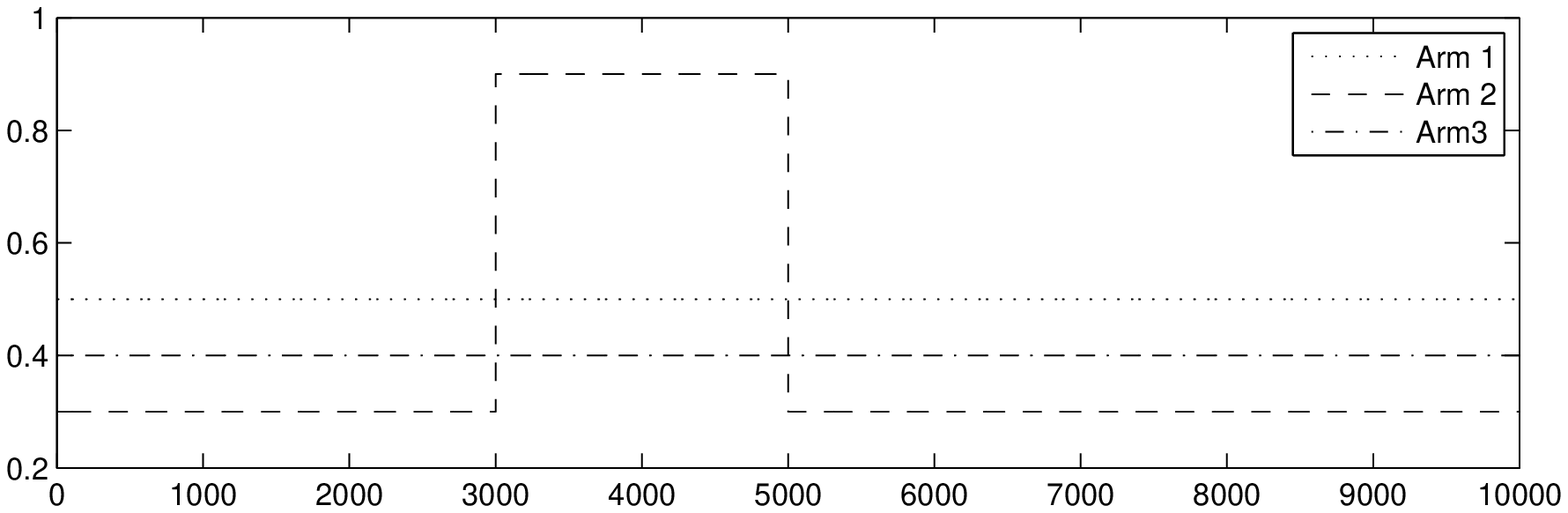} &  \includegraphics[width=7cm,height=2cm]{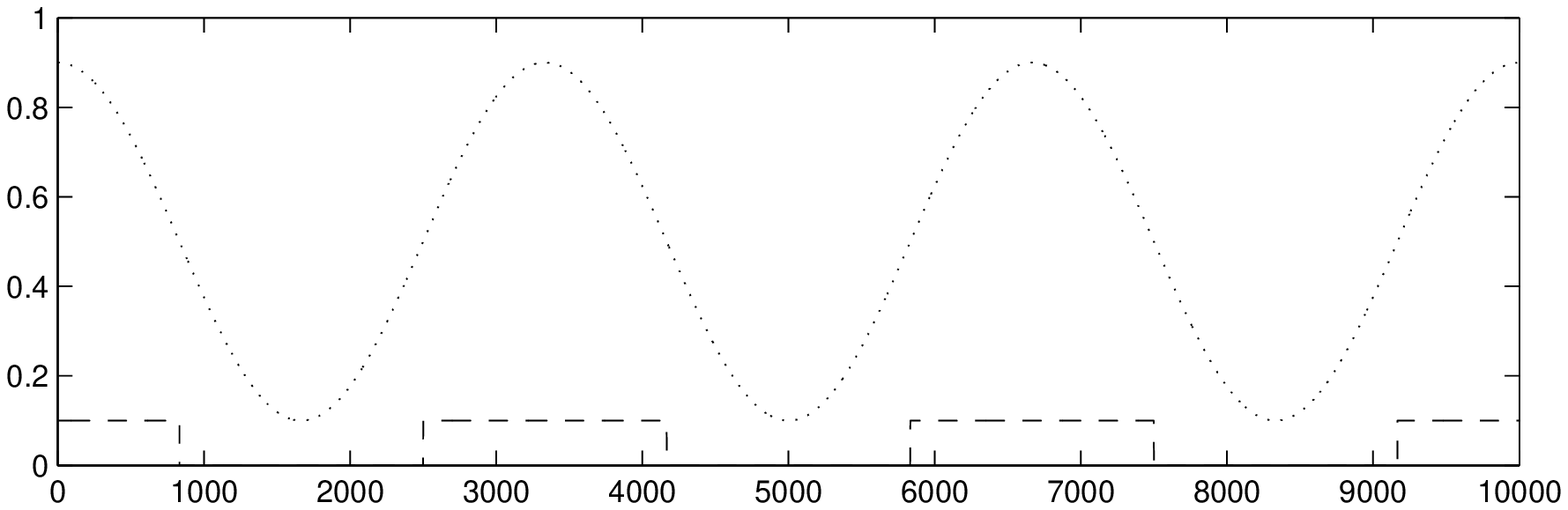}\\
 \includegraphics[width=7cm]{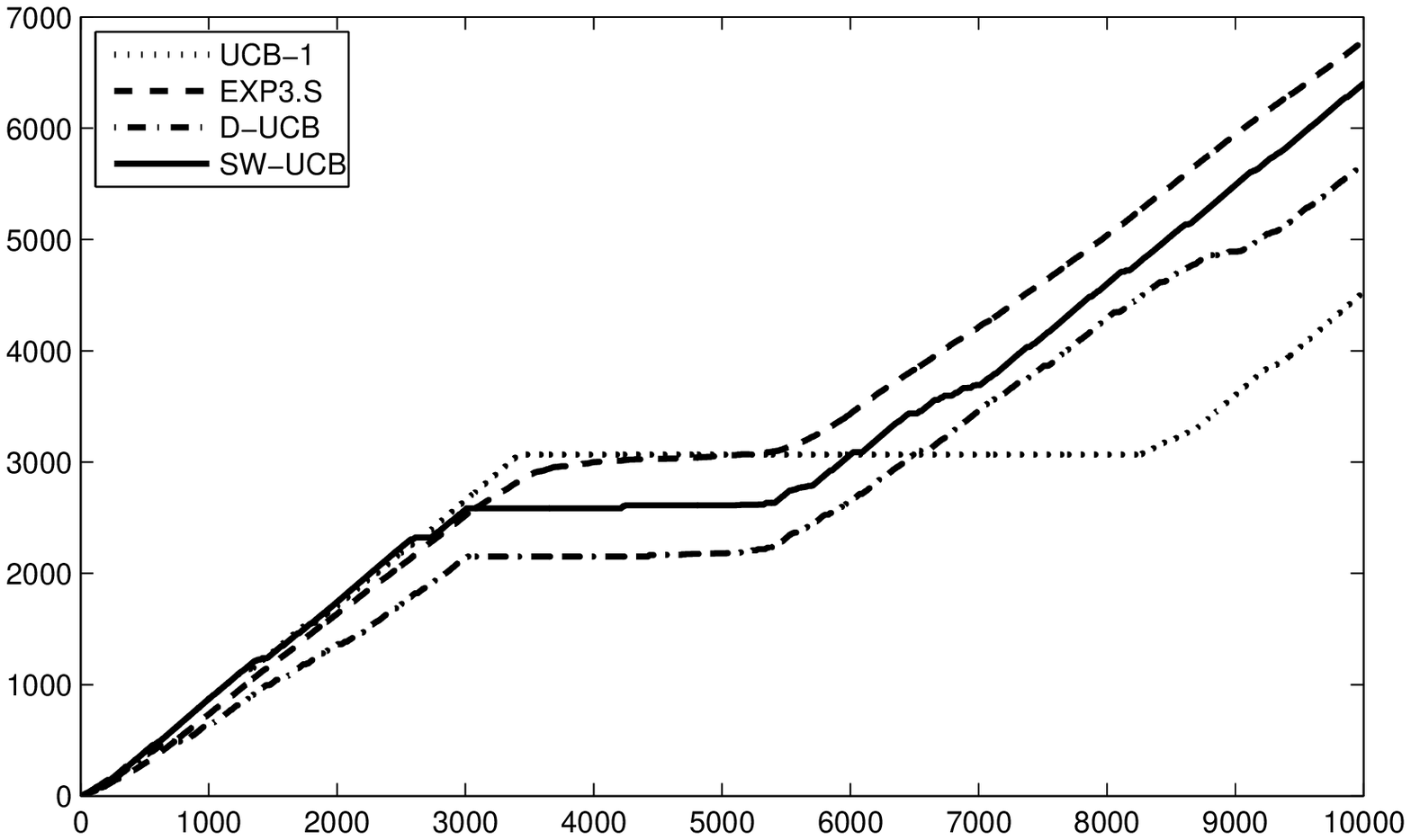} &  \includegraphics[width=7cm]{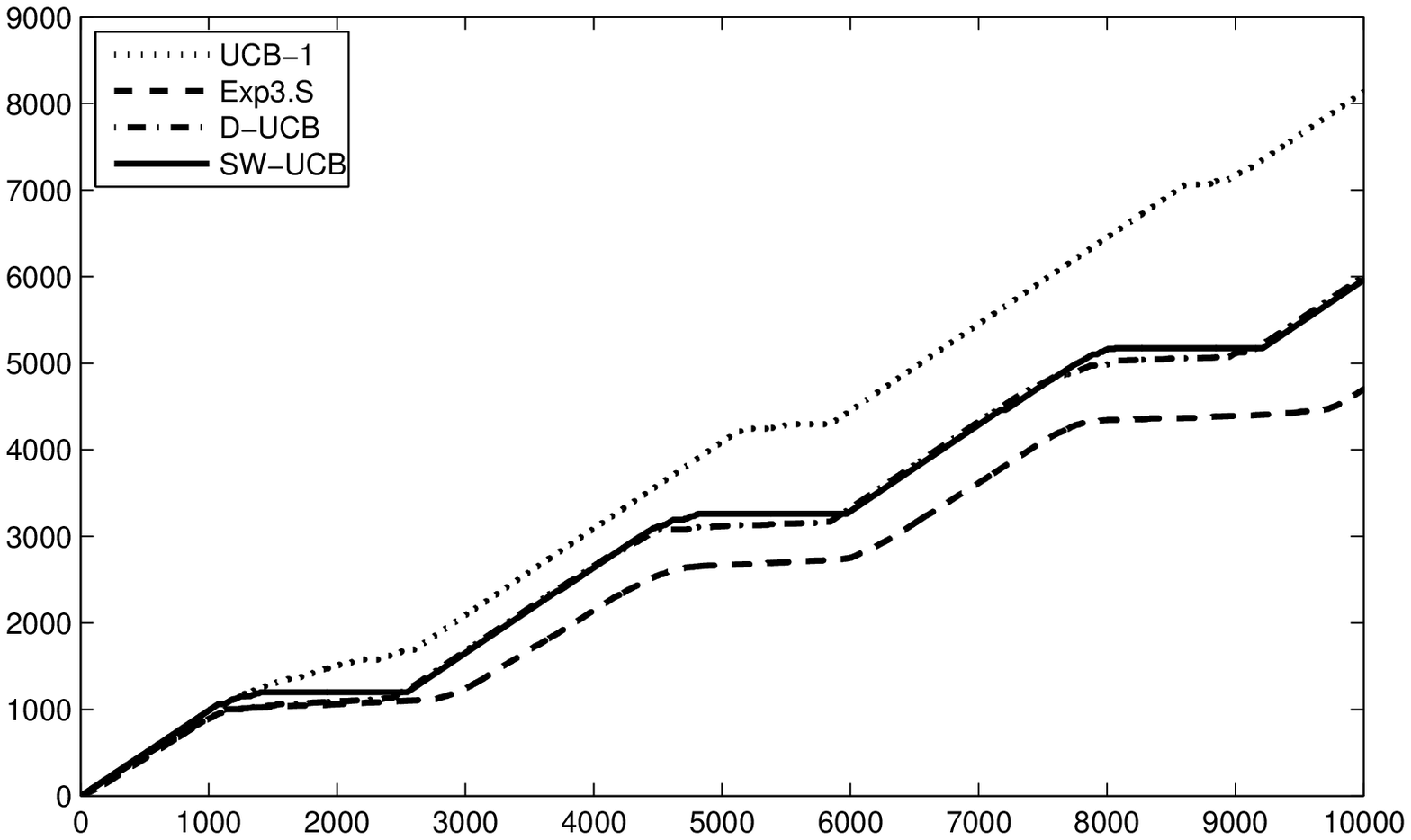}\\
 \includegraphics[width=7cm]{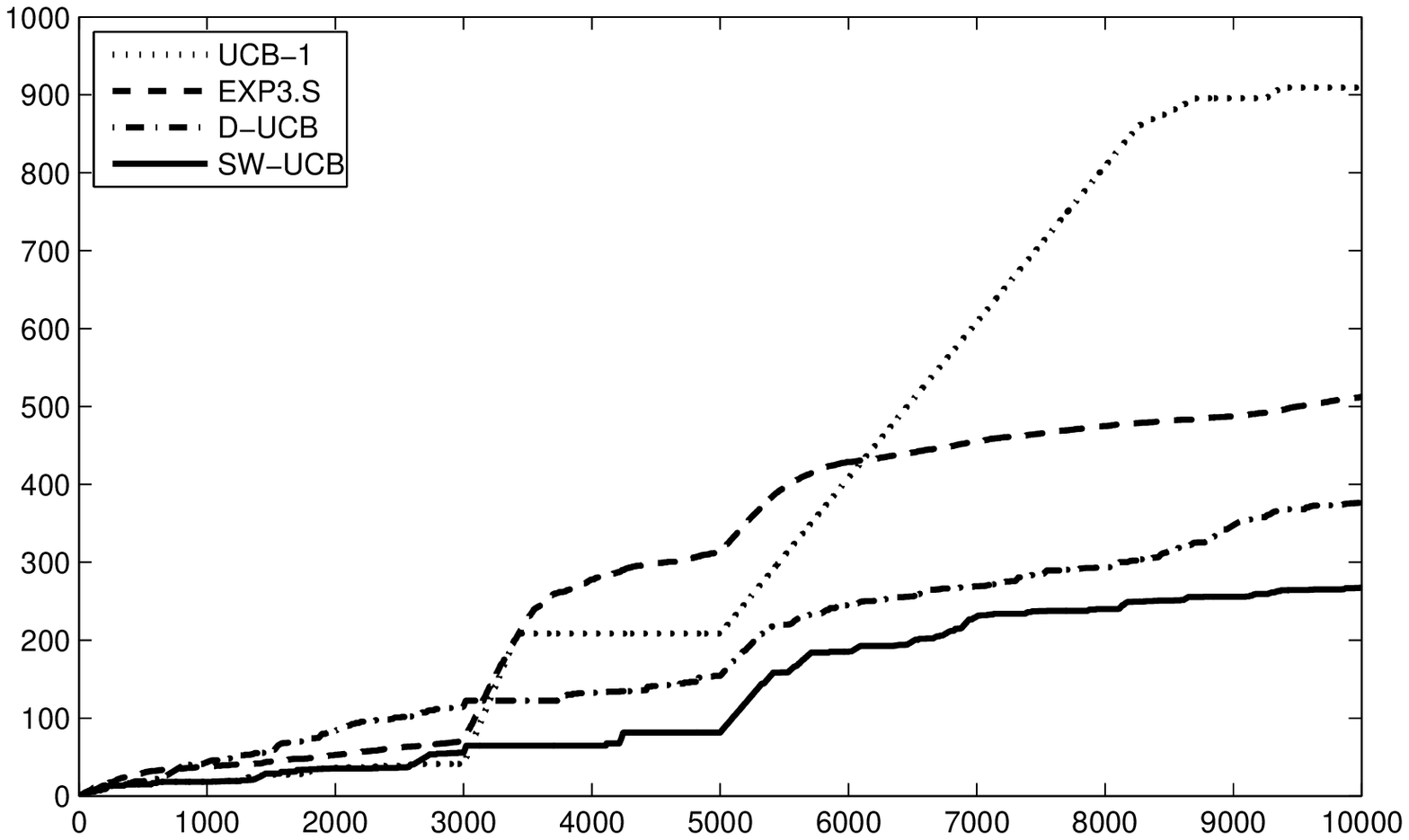} & \includegraphics[width=7cm]{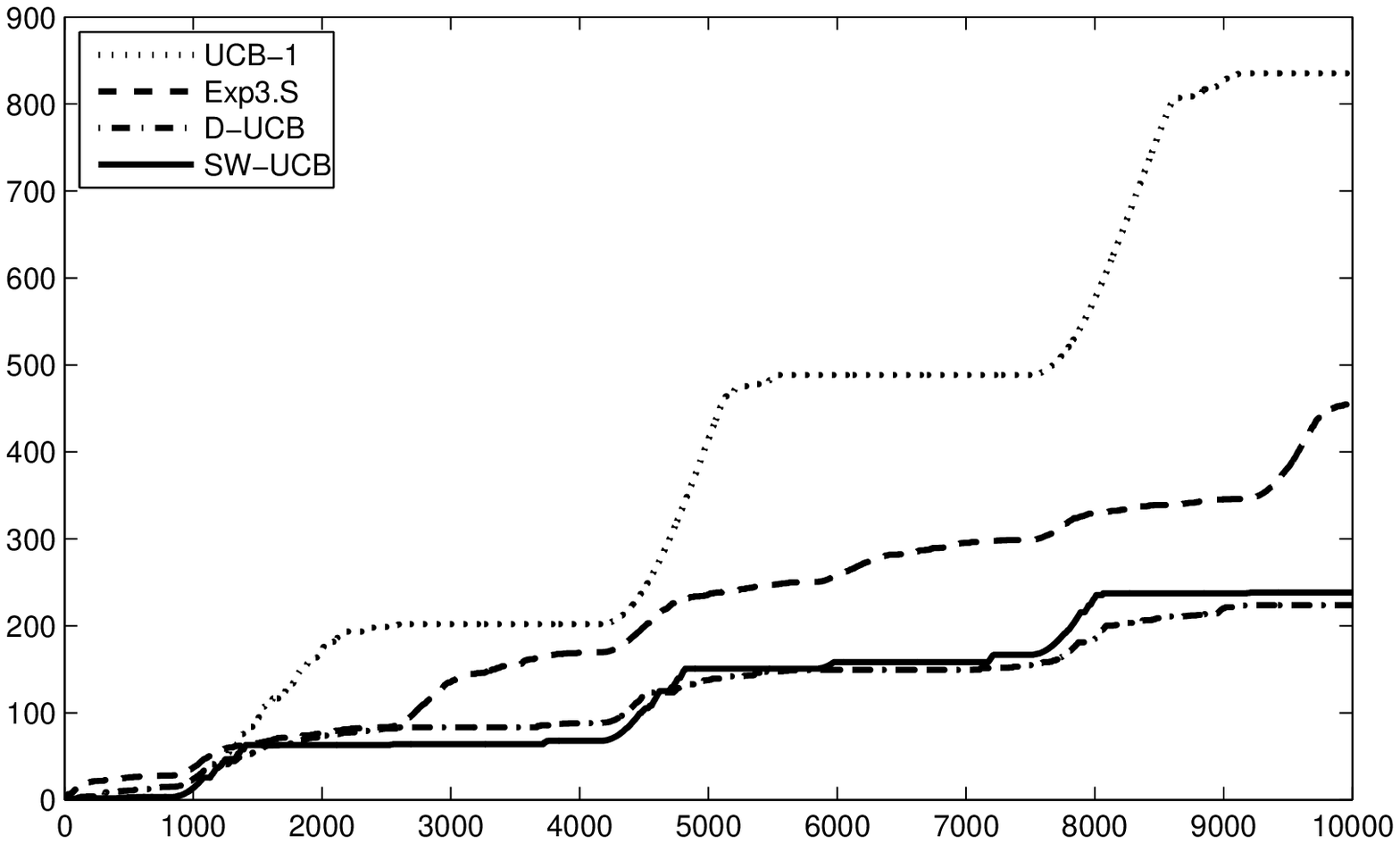}\\
\end{tabular}
 \caption{
Left panel:
 Bernoulli MAB problem with two swaps.
 Upper: evolution of the probability of having a reward 1 for each arm;
 Middle: cumulative frequency of arm 1 pulls for each policy.
 Below: cumulative regret of each policy.
Right panel:
Bernoulli MAB problem  with periodic rewards:
Upper: evolution of the probability of having a reward 1 for arm 1 (and time intervals when it should be played);
Middle: cumulative frequency of arm 1 pulls for each policy.
Below: cumulative regret of each policy.
}
\label{fig1}
\end{figure}

We consider here two settings.
In the first example, there are $K=3$ arms and the time horizon is set to $T=10^4$.
The agent goal is to minimize the expected regret. The rewards of arm $i \in \{1,\dots,K\}$ at time $t$ are independent Bernoulli random variables with success probability $p_t(i)$, with $p_t(1)=0.5$, $p_t(2)=0.3$ and for $t\in\{1,\dots, T\}$: $$p_t(3)=\left\{\begin{array}{rcl}0.4 & \hbox{ for }& t<3000 \hbox{ or } t\geq 5000,\\0.9 & \hbox{ for }&  3000\leq t <5000.\end{array}\right.$$
As one may notice, the optimal policy for this bandit task is to select arm 1 before the first breakpoint ($t=3000$) and after the second breakpoint ($t=5000$). In the left panel of Figure \ref{fig1}, we represent the evolution of two criteria in function of $t$: the number of times policy 1 has been played, and the cumulated regret (bottom plot).
These two measures are obviously related, but they are not completely equivalent as sub-optimal arms can yield relatively high rewards.
We compare the UCB-1 algorithm with $\xi=\frac{1}{2}$,
the EXP3.S algorithm described in \cite{auer:cesabianchi:freund:schapire:2002} with the tuned parameters given in Corollary 8.3 (with the notations of this paper $\alpha= T^{-1}$ and $\gamma= \sqrt{K(\Upsilon_T \log(KT)+\rme)/[(\rme-1)T]}$ with $\Upsilon_T=2$),
the D-UCB algorithm with $\xi=1/2$ and $\gamma=1-1/4\sqrt{T}$ and the SW-UCB with $\xi=1/2$ and $\tau=4\sqrt{n\log T}$.
The parameters are tuned to obtain roughly optimal performance for the chosen horizon $T$ and the number of breakpoints.

As can be seen in Figure \ref{fig1} (and as consistently observed over the simulations), D-UCB performs almost as well as SW-UCB. Both of them waste significantly less time than EXP3.S and UCB-1 to detect the breakpoints, and quickly concentrate their pulls on the optimal arm.
Observe that policy UCB-1, initially the best, reacts very fast to the first breakpoint ($t=3000$), as the confidence interval for arm $3$ at this step is very loose. On the contrary, it takes a very long time after the second breakpoint ($t=5000$) for UCB-1 to play arm $1$ again.

In the second example, there are $K=2$ arms, the rewards are still Bernoulli random variables with parameters $p_t(i)$ but are in persistent, continuous evolution.
Arm 2 is taken as a reference ($p_t(2)=1/2$ for all $t$), and the parameter of arm 1 evolves periodically as: $p_t(1)=0.5+0.4\cos\left( 6\pi Rt /T \right)$.
Hence, the best arm to pull evolves cyclically and the transitions are smooth (regularly, the two arms are statistically indistinguishable).
The middle plot in the right panel of Figure \ref{fig1} represents the cumulative frequency of arm 1 pulls: D-UCB, SW-UCB and, to a lesser extent, EXP3.S   track  the cycles, while UCB-1  fails to identify the best current arm. 
Below, the evolutions of the cumulative regrets under the four policies are shown: in this continuously evolving environment, the performance of 
D-UCB and SW-UCB are almost equivalent while UCB-1 and the Exp3.S algorithms accumulate larger regrets.

\section{Conclusion and perspectives}
This paper theoretically establishes that the UCB policies can also be successfully adapted to cope with non-stationary environments.  
The upper bound of the SW-UCB in abruptly changing environment matches the upper bounds of the Exp3.S algorithm (i.e. $O(\sqrt{T \log(T)})$), showing that UCB policies can be at least as good as the softmax methods. In practice, numerical experiments also support this finding. For the two examples considered in this paper, the D-UCB and SW-UCB policies outperform the optimality tuned version of the Exp3.S algorithm.

The focus of this paper is on abruptly changing environment, but it is believed that the theoretical tools developed to handle the non-stationarity can be applied in different contexts. In particular, using a similar  bias-variance decomposition of the discounted or windowed-rewards, the analysis of continuously evolving reward distributions can be done (and will be reported in a forthcoming paper). Furthermore, Theorems \ref{th:concentrationR}
and \ref{th:concentrationRmaxi}, dealing with concentration inequality for discounted martingale transforms, are powerful tools of independent interest. 

As the previously reported Exp3.S algorithm, the performance of the proposed policy depends on tuning parameters, the discount factor for D-UCB and the window size for SW-UCB. These tuning parameters may be adaptively set, using data-driven approaches, such as the one proposed in \cite{hartland:gelly:baskiotis:teytaud:sebag:2006}. This is the subject of on-going research.

\appendix

\section{A Hoeffding-type inequality for self-normalized means with a random number of summands}
\label{sec:self-normalizedhoeffing}
Let $(X_t)_{t\geq 1}$  be a sequence of non-negative independent bounded random variables defined on a probability space $(\Omega, \A, \P)$. We denote by $B$ the upper bound, $X_t \in [0,B]$, $\P$-a.s. and by $\mu_t$ its expectation $\mu_t= \E[X_t]$.
Let $\F_t$ be an increasing sequence of $\sigma$-fields of $\A$ such that for each $t$, $\sigma(X_1\ldots, X_t)\subset \F_t$ and for $s>t$, $X_s$ is independent from $\F_t$.
Consider a previsible sequence $(\epsilon_t)_{t\geq 1}$ of Bernoulli variables (for all $t > 0$, $\epsilon_t$ is $\F_{t-1}$-measurable).
Denote by $\phi_t$ the Cramer transform of $X_t$: for $\lambda\in\Rset$,
$$\phi_t(\lambda)=\log \E[\exp(\lambda X_t)] \eqsp.$$
For $\gamma \in [0,1)$
, consider the following random variables
\begin{equation}
\Reward{t}{\gamma}  =\sum_{s=1}^t \gamma^{t-s}X_s\epsilon_s \eqsp , \qquad
\EReward{t}{\gamma} = \sum_{s=1}^t \gamma^{t-s}\mu_s\epsilon_s\eqsp,\qquad
\Nplay{t}{\gamma} = \sum_{s=1}^t \gamma^{t-s}\epsilon_s\eqsp.
\end{equation}
Let also
$$\Nplay[sum][]{t}{\gamma}=\sum_{s=1}^t \gamma^{t-s}=\left\{\begin{array}{ccc}\frac{1-\gamma^t}{1-\gamma}& \hbox{if}&\gamma<1,\\ t&\hbox{if}&\gamma=1.\end{array}\right.$$
\begin{theorem}
For all integers $t$ and all $\delta>0$,
\label{th:concentrationR}
\begin{equation*}
 \P\left(\frac{\Reward{t}{\gamma}-\EReward{t}{\gamma}}{\sqrt{\Nplay{t}{\gamma^2}}} > \delta\right)\\\leq
\left\lceil\frac{\log \Nplay[sum][]{t}{\gamma}}{\log (1+\eta)}\right\rceil \exp\left(-\frac{2\delta^2}{B^2}\left(1-\frac{\eta^2}{16}\right)\right)
\end{equation*}
for all $\eta>0$.
\end{theorem}
\begin{remark}
Actually, we prove the slightly stronger inequality:
\begin{equation}
\label{eq:concentrationRoptimale}
\P\left(\frac{\Reward{t}{\gamma}-\EReward{t}{\gamma}}{\sqrt{\Nplay{t}{\gamma^2}}} > \delta\right)
\leq \left\lceil\frac{\log \Nplay[sum][]{t}{\gamma}}{\log (1+\eta)}\right\rceil \\\times \exp\left(-\frac{8\delta^2}{B^2\left( (1+\eta)^{1/4}+ (1+\eta)^{-1/4} \right)^2}\right)\eqsp.
\end{equation}
\end{remark}
\begin{proof}
First observe that we can assume $\epsilon_t=1$, since otherwise
$(\Reward{t}{\gamma}-\EReward{t}{\gamma})/\sqrt{\Nplay{t}{\gamma^2}} = (\Reward{t-1}{\gamma}-\EReward{t-1}{\gamma})/\sqrt{\Nplay{t-1}{\gamma^2}}$ and the result follows from a simple induction.
Second, note that for every positive $\lambda$ and for every  $u<t$, since $\epsilon_{u+1}$ is predictible, and since $X_{u+1}$ is independent from $\F_{u}$,
$$
\CPE{\exp\left( \lambda X_{u+1}\epsilon_{u+1}\right)}{\F_{u}}
= \exp\left(\phi_{u+1}\left(\lambda  \epsilon_{u+1}\right)\right)
= \exp\left(\phi_{u+1}\left(\lambda \right) \epsilon_{u+1}\right)\eqsp.
$$
Hence, as $\Reward{u+1}{\gamma} = \gamma \Reward{u}{\gamma}+X_{u+1}\epsilon_{u+1}$,
$$\E\left[ \exp\left( \lambda\Reward{u+1}{\gamma}-\sum_{s=1}^{u+1}\phi_s\left( \lambda\gamma^{u+1-s} \right)\epsilon_s \right)\right]=\E\left[ \exp\left( \lambda\gamma\Reward{u}{\gamma}-\sum_{s=1}^{u}\phi_s\left( (\lambda\gamma)\gamma^{u-s} \right)\epsilon_s \right)\right]\eqsp.
$$
As $\phi(0)=0$, this proves by induction that
$$\E\left[\exp\left( \lambda\Reward{t}{\gamma}-\sum_{s=1}^{t}\phi_s\left( \lambda\gamma^{t-s} \right)\epsilon_s\right)\right]=1\eqsp.$$
It is easily verified (see e.g. \cite[Lemma 8.1]{MR1383093}) that under the stated assumptions, for all positive $\lambda$,
\begin{equation}
\label{eq:BoundHoeffding}
\phi_s(\lambda) \leq \lambda \mu_s+ B^2 \lambda^2/8 \eqsp,
\end{equation}
showing that
\begin{equation*}
\E\left[\exp\left( \lambda\left\{ \Reward{t}{\gamma} - \EReward{t}{\gamma} \right\} - \frac{B^2}{8}\lambda^2 \Nplay{t}{\gamma^2} \right)\right] \leq 1.
\end{equation*}
Hence, for any $x>0$, the Markov inequality yields
\begin{multline*}
\P\Bigg(\frac{\Reward{t}{\gamma}-\EReward{t}{\gamma}}{\sqrt{\Nplay{t}{\gamma^2}}} > \frac{x}{\lambda\sqrt{\Nplay{t}{\gamma^2}}} + \frac{\lambda B^2\sqrt{\Nplay{t}{\gamma^2}}}{8} \Bigg)
\\ =\P\left(\exp\left( \lambda\left( \Reward{t}{\gamma} - \EReward{t}{\gamma} \right) - \frac{B^2}{8}\lambda^2 \Nplay{t}{\gamma^2} \right)\geq \rme^x \right)
 \leq \exp(-x)\eqsp.
\end{multline*}
Now, take $\eta>0$, let $D=\left\lceil\frac{\log \Nplay[sum][]{t}{\gamma}}{\log (1+\eta)}\right\rceil$ and, for every integer $k\in\{1,\dots,D\}$, define
$$\lambda_k=\sqrt{\frac{8x}{B^2(1+\eta)^{k-\frac{1}{2}}}}.$$
Elementary algebra shows that for all $z$ such that $(1+\eta)^{k-1}\leq z \leq (1+\eta)^{k}$, we have
\begin{equation}
\sqrt{\frac{(1+\eta)^{k-\frac{1}{2}}}{z}} + \sqrt{\frac{z}{(1+\eta)^{k-\frac{1}{2}}}} \leq (1+\eta)^{1/4}+ (1+\eta)^{-1/4}
\label{eq:elemalgebra}\end{equation}
Thus, if $(1+\eta)^{k-1}\leq\Nplay{t}{\gamma^2}\leq(1+\eta)^{k}$, then
\begin{align*}
\frac{x}{\lambda_k\sqrt{\Nplay{t}{\gamma^2}}} + \frac{B^2}{8}\lambda_k\sqrt{\Nplay{t}{\gamma^2}}&= B\sqrt{\frac{x}{8}}\left( \sqrt{\frac{(1+\eta)^{k-\frac{1}{2}}}{\Nplay{t}{\gamma^2}}} + \sqrt{\frac{ \Nplay{t}{\gamma^2}}{(1+\eta)^{k-\frac{1}{2}}}} \right)\\
&\leq B\sqrt{\frac{x}{8}}\left( (1+\eta)^{1/4}+ (1+\eta)^{-1/4} \right).
\end{align*}
Therefore, as $\epsilon_t=1$ we have $1\leq N_t(\gamma^2)\leq (1+\eta)^D$ and
\begin{multline*}
\left\{\frac{\Reward{t}{\gamma}-\EReward{t}{\gamma}}{\sqrt{\Nplay{t}{\gamma^2}}} > B\sqrt{\frac{x}{8}}\left( (1+\eta)^{1/4}+ (1+\eta)^{-1/4} \right)\right\}\\
\subset \bigcup_{k=1}^D\left\{\frac{\Reward{t}{\gamma}-\EReward{t}{\gamma}}{\sqrt{\Nplay{t}{\gamma^2}}} > \frac{x}{\lambda_k\sqrt{\Nplay{t}{\gamma^2}}} + \frac{\lambda_k B^2\sqrt{\Nplay{t}{\gamma^2}}}{8}\right\}\eqsp .
\end{multline*}
The union bound thus implies that:
\begin{multline*}
\P\left( \frac{\Reward{t}{\gamma}-\EReward{t}{\gamma}}{\sqrt{\Nplay{t}{\gamma^2}}} > B\sqrt{\frac{x}{8}}\left( (1+\eta)^{1/4}+ (1+\eta)^{-1/4} \right) \right)\\
\leq  \sum_{k=1}^D\P\left( \frac{\Reward{t}{\gamma}-\EReward{t}{\gamma}}{\sqrt{\Nplay{t}{\gamma^2}}} > \frac{x}{\lambda_k\sqrt{\Nplay{t}{\gamma^2}}} + \frac{\lambda_k B^2\sqrt{\Nplay{t}{\gamma^2}}}{8}\right)
\leq D\exp(-x)\eqsp .
\end{multline*}
For $\delta=B\sqrt{\frac{x}{8}}\left( (1+\eta)^{1/4}+ (1+\eta)^{-1/4} \right)$, this yields
\begin{equation*}
\P\left(\frac{\Reward{t}{\gamma}-\EReward{t}{\gamma}}{\sqrt{\Nplay{t}{\gamma^2}}} > \delta\right)\\
\leq D \exp\left(-\frac{8\delta^2}{B^2\left( (1+\eta)^{1/4}+ (1+\eta)^{-1/4} \right)^2}\right)\eqsp.
\end{equation*}
The conclusion follows, as it is easy to see that, for all $\eta>0$,
\begin{equation}
\frac{4}{\left( (1+\eta)^{1/4}+ (1+\eta)^{-1/4} \right)^2}\geq 1-\frac{\eta^2}{16} \eqsp.
\label{eq:majeta}
\end{equation}
\end{proof}

\begin{remark}
For example, taking $\eta=0.3$ in (\ref{eq:concentrationRoptimale}) yields
\begin{equation*}
 \P\left(\frac{\Reward{t}{\gamma}-\EReward{t}{\gamma}}{\sqrt{\Nplay{t}{\gamma^2}}} > \delta\right)
\leq
4\log\Nplay[sum][]{t}{\gamma} \rme^{-\frac{1.99\delta^2}{B^2}} \eqsp.
\end{equation*}
Classical Hoeffding bounds for deterministic $\epsilon_s$ yield an upper-bound in
$$\P\left( \frac{\Reward{t}{\gamma}-\EReward{t}{\gamma}}{\sqrt{\Nplay{t}{\gamma^2}}} > \delta\right) \leq \rme^{-\frac{2\delta^2}{B^2}}$$
for all positive $t$.
The factor behind the exponential and the very slightly larger exponent are the price to pay for the presence of random $\epsilon_s$. Theorem \ref{th:concentrationR} is maybe sub-optimal, but it is possible to show that for all $\delta>0$ and for an appropriate choice of the previsible sequence $(\epsilon_s)_{s\geq 1}$
$$\P\left( \frac{\Reward{t}{\gamma}-\EReward{t}{\gamma}}{\sqrt{\Nplay{t}{\gamma^2}}} > \delta\right) \to 1$$
as $t$ goes to infinity.

\end{remark}

If all variables $X_t$ have the same expectation $\mu$, taking $\gamma=1$ in Theorem \ref{th:concentrationR} immediately leads to  the following corollary:
\begin{corollary}
\label{cor:concentrationMR}
For all integers $t$ and $\tau$,
\begin{equation*}\P\left(\frac{\sum_{s=(t-\tau+1)\wedge 1}^t (X_s -\mu)\epsilon_s}{\sqrt{\sum_{s=(t-\tau+1)\wedge 1}^t\epsilon_s}}> \delta\right)\\\leq
 \left\lceil\frac{\log (t\wedge \tau)}{\log(1+\eta)}\right\rceil \exp\left(-\frac{2\delta^2}{B^2}\left(1-\frac{\eta^2}{16}\right)\right)
 \end{equation*}
\end{corollary}

\section{A maximal inequality for self-normalized means with a random number of summands}
\label{sec:self-normalizedmaxihoeffing}
In this section, we prove a stronger version of Theorem \ref{th:concentrationR}: we upper-bound the probability that, at some time $t$, the average reward deviates from its expectation. We keep the same notations as in Section \ref{sec:self-normalizedhoeffing}.
\begin{theorem}
For all positive integer $T$ and all $\delta>0$,
\label{th:concentrationRmaxi}
\begin{equation*}
 \P\left( \sup_{1\leq t\leq T}\frac{\Reward{t}{\gamma}-\EReward{t}{\gamma}}{\sqrt{\Nplay{t}{\gamma^2}}} > \delta\right)\\\leq
\left\lceil \frac{\log \left( \gamma^{-2T}\Nplay[sum][]{T}{\gamma^2} \right)}{\log(1+\eta)} \right\rceil \exp\left(-\frac{2\delta^2}{B^2}\left(1-\frac{\eta^2}{16} \right)\right)\eqsp.
\end{equation*}
for all $\eta>0$.
\end{theorem}
\begin{remark}
Note that if $\gamma<1$, then
$$\log \left(\gamma^{-2T}\Nplay[sum][]{T}{\gamma^2} \right)\leq \frac{2T(1-\gamma)}{\gamma} + \log\frac{1}{1-\gamma^2}\eqsp,$$
while for $\gamma=1$ we have:
$$\log \left(\gamma^{-2T}\Nplay[sum][]{T}{\gamma^2} \right)=\log T.$$
\end{remark}
\begin{remark}
Classical Hoeffding bounds for deterministic $\epsilon_s$ yield an upper-bound in
$$\P\left( \frac{\Reward{t}{\gamma}-\EReward{t}{\gamma}}{\sqrt{\Nplay{t}{\gamma^2}}} > \delta\right) \leq \exp(-2\delta^2)$$
for all positive $t$.
The factor behind the exponential (depending on $T$ and $\epsilon$) and the very slightly larger exponant are the price to pay for uniformity in $t$.
For example, taking $\eta=0.3$ yields
\begin{equation*}
 \P\left( \sup_{1\leq t\leq T}\frac{\Reward{t}{\gamma}-\EReward{t}{\gamma}}{\sqrt{\Nplay{t}{\gamma^2}}} > \delta\right)\\\leq
\left\lceil 4\log \left(\gamma^{-2T}\Nplay[sum][]{T}{\gamma^2}\right) \right\rceil \exp\left(-\frac{1.99\delta^2}{B^2}\right)\eqsp.
\end{equation*}
\end{remark}

\begin{proof}
For $\lambda>0$, define
\begin{equation}
Z^\lambda_t  =  \exp\left( \lambda\gamma^{-t} \Reward{t}{\gamma} - \sum_{s=1}^{t}\phi_s\left( \lambda \gamma^{-s}\right)\epsilon_s \right) \eqsp.
\end{equation}
Note that
\[
\CPE{\exp\left( \lambda\gamma^{-t} X_t\epsilon_t\right)}{\F_{t-1}}
= \exp\left( \epsilon_t\phi_t\left(\lambda \gamma^{-t}\right)\right).
\]
Since $\gamma^{-t} \Reward{t}{\gamma}= \gamma^{-(t-1)} \Reward{t-1}{\gamma} + \gamma^{-t} X_t \epsilon_t$,
we may therefore write
\begin{equation*}
\CPE{\exp(\lambda \gamma^{-t} \Reward{t}{\gamma})}{\F_{t-1}}\\
= \exp\left( \lambda \gamma^{-(t-1)} \Reward{t-1}{\gamma} \right) \exp\left(\epsilon_t \phi_t(\lambda \gamma^{-t})\right) \eqsp,
\end{equation*}
showing that $\{Z^\lambda_t \}$ is a martingale adapted to the filtration $\F= \{ \mathcal{F}_t, t \geq 0 \}$.
As already mentionned (see e.g. \cite[Lemma 8.1]{MR1383093}), under the stated assumptions
\begin{equation*}
\phi_t(\lambda) \leq \lambda \mu_t+ B^2 \lambda^2/8 \eqsp,
\end{equation*}
showing that for all $\lambda>0$,
\begin{equation}
W^\lambda_t  =  \exp\left( \lambda\gamma^{-t} \Reward{t}{\gamma} - \right.\\\left. \lambda\gamma^{-t}\EReward{t}{\gamma} - (B^2/8)\lambda^2 \gamma^{-2t}\Nplay{t}{\gamma^2} \right)
\end{equation}
is a super-martingale. Hence, for any $x>0$ we have
\begin{equation}
\label{eq:majprobaWlambda}
\P\left(\sup_{
1\leq t\leq T}W^\lambda_t\geq \exp(x) \right) \leq \exp(-x)\eqsp.
\end{equation}
On the other hand, note that
\begin{equation}
\left\{W^\lambda_t > \exp(x)\right\}
=\left\{ \frac{\Reward{t}{\gamma}-\EReward{t}{\gamma}}{\sqrt{\Nplay{t}{\gamma^2}}}>\frac{x\gamma^t}{\lambda\sqrt{\Nplay{t}{\gamma^2}}} + \frac{B^2}{8}\lambda\gamma^{-t}\sqrt{\Nplay{t}{\gamma^2}}\right\}\eqsp.\label{eq:Wgdimpliesecartgrand}
\end{equation}
Now, let $D=\left\lceil\frac{\log \left(\gamma^{-2T}\Nplay[sum][]{T}{\gamma^2} \right)}{\log (1+\eta)}\right\rceil$ and for every integer $k\in\{1,\dots,D\}$, define
$$\lambda_k=\sqrt{\frac{8x}{B^2(1+\eta)^{k-\frac{1}{2}}}}.$$
Thus, if $(1+\eta)^{k-1}\leq \gamma^{-2t}\Nplay{t}{\gamma^2}\leq(1+\eta)^{k}$, then using Equation  \eqref{eq:elemalgebra} yields:
\begin{align*}
\frac{x\gamma^t}{\lambda_k\sqrt{\Nplay{t}{\gamma^2}}} + \frac{B^2}{8}\lambda_k\gamma^{-t}\sqrt{\Nplay{t}{\gamma^2}}
&= B\sqrt{\frac{x}{8}}\left( \sqrt{\frac{(1+\eta)^{k-\frac{1}{2}}}{ \gamma^{-2t}\Nplay{t}{\gamma^2}}} + \sqrt{\frac{ \gamma^{-2t}\Nplay{t}{\gamma^2}}{(1+\eta)^{k-\frac{1}{2}}}} \right)\\
&\leq B\sqrt{\frac{x}{8}}\left( (1+\eta)^{1/4}+ (1+\eta)^{-1/4} \right),
\end{align*}
which proves, using Equation (\ref{eq:Wgdimpliesecartgrand}), that
\begin{multline*}
\left\{\frac{\Reward{t}{\gamma}-\EReward{t}{\gamma}}{\sqrt{\Nplay{t}{\gamma^2}}} > B\sqrt{\frac{x}{8}}\left( (1+\eta)^{1/4}+ (1+\eta)^{-1/4} \right)\right\} \\\subset \bigcup_{k=1}^D
\left\{\frac{\Reward{t}{\gamma}-\EReward{t}{\gamma}}{\sqrt{\Nplay{t}{\gamma^2}}} > \frac{x\gamma^t}{\lambda_k\sqrt{\Nplay{t}{\gamma^2}}} + \frac{B^2}{8}\lambda_k\gamma^{-t}\sqrt{\Nplay{t}{\gamma^2}}\right\} \subset \left\{W^{\lambda_k}_t > \exp(x)\right\}\eqsp .
\end{multline*}
But as $$\sup_{1\leq t\leq T}\frac{\Reward{t}{\gamma}-\EReward{t}{\gamma}}{\sqrt{\Nplay{t}{\gamma^2}}} = \sup_{1\leq t\leq T, \epsilon_t=1}\frac{\Reward{t}{\gamma}-\EReward{t}{\gamma}}{\sqrt{\Nplay{t}{\gamma^2}}}\eqsp,$$
we assume $\epsilon_t=1$ and thus $1\leq \Nplay{t}{\gamma^2} \leq (1+\eta)^D$.
Hence, thanks to Equation (\ref{eq:majprobaWlambda}) we obtain:
\begin{multline*}
\P\Bigg( \bigcup_{1\leq t\leq T} \Bigg\{\frac{\Reward{t}{\gamma}-\EReward{t}{\gamma}}{\sqrt{\Nplay{t}{\gamma^2}}}
> B\sqrt{\frac{x}{8}}\left( (1+\eta)^{1/4}+ (1+\eta)^{-1/4} \right) \Bigg\} \Bigg)\\
\leq \P \left( \bigcup_{1\leq t\leq T}\bigcup_{1\leq k\leq D}\Big\{ W^{\lambda_k}_t > \exp(x)\Big\}\right)
= \P \left( \bigcup_{1\leq k\leq D}\bigcup_{1\leq t\leq T}\Big\{ W^{\lambda_k}_t > \exp(x)\Big\} \right)
\leq D \exp(-x)\eqsp.
\end{multline*}
For $$\delta=B\sqrt{\frac{x}{8}}\left( (1+\eta)^{1/4}+ (1+\eta)^{-1/4} \right),$$
and using Equation \eqref{eq:majeta}, this yields
\begin{align*}
\P\left(\frac{\Reward{t}{\gamma}-\EReward{t}{\gamma}}{\sqrt{\Nplay{t}{\gamma^2}}} > \delta\right)&
\leq D \exp\left(-\frac{8\delta^2}{B^2\left( (1+\eta)^{1/4}+ (1+\eta)^{-1/4} \right)^2}\right)\\
&\leq D\exp\left(-\frac{2\delta^2}{B^2}\left(1-\frac{\eta^2}{16} \right)\right)\eqsp.
\end{align*}
 \end{proof}

\section{Technical results}
\label{sec:technicalresults}
\begin{lemma}
For any $i\in\{1,\dots,K\}$ and for any positive integer $\tau$, let $\Nplay[][i]{t-\tau:t}{1}= \sum_{s=t-\tau+1}^t \1{I_t=i}$.
Then for any positive $m$,
 \label{lem:majabadondiscount}
$$
\sum_{t=K+1}^T \1{I_t=i, \Nplay[][i]{t-\tau:t}{1}< m} \leq K\lceil T/\tau\rceil m \eqsp.
$$
\end{lemma}

\begin{proof}
\begin{equation*}
\sum_{t=1}^T \1{I_t=i, \Nplay[][i]{t-\tau:t}{1}< m} \\\leq \sum_{j=1}^{\lceil T/\tau\rceil} \sum_{t=(j-1)\tau+1}^{j\tau} \1{I_t=i, \Nplay[][i]{t-\tau:t}{1}< m}.
\end{equation*}
For any given $j\in\{1,\ldots,\lceil T/\tau\rceil\}$, either
$\sum_{t=(j-1)\tau+1}^{j\tau} \1{I_t=i, \Nplay[][i]{t-\tau:t}{1}< m}=0$
or there exists an index $t \in \{(j-1)\tau+1,\dots, j\tau\}$ such that $I_t=i$, $\Nplay[][i]{t-\tau:t}{1}< m$.
In this case, we put $t_j=\max\{t\in\{(j-1)\tau+1,\dots, j\tau\}:
I_t=i, \Nplay[][i]{t-\tau:t}{1}< m\}$, the last time this condition is met in the $j$-th block.
Then,
\begin{multline*}
\sum_{t=(j-1)\tau+1}^{j\tau} \1{I_t=i, \Nplay[][i]{t-\tau:t}{1}< m}
= \sum_{t=(j-1)\tau+1}^{t_j} \1{I_t=i, \Nplay[][i]{t-\tau:t}{1}< m} \\
\leq  \sum_{t=t_j-\tau+1}^{t_j} \1{I_t=i, \Nplay[][i]{t-\tau:t}{1}< m}
\leq  \sum_{t=t_j-\tau+1}^{t_j} \1{I_t=i}
= \Nplay[][i]{t_j-\tau:t_j}{1}
< m.
\end{multline*}
\end{proof}

\begin{corollary}
 \label{cor:majabadondiscount}
For any $i\in\{1,\dots,K\}$, any integers $\tau \geq 1$ and $A > 0$,
$$
\sum_{t=K+1}^T \1{I_t=i, \Nplay[][i]{t}{\gamma}< A} \leq K\lceil T/\tau\rceil A \gamma^{-\tau} \eqsp.
$$
\end{corollary}

\begin{proof}
Simply note that
\begin{equation}
\sum_{t=K+1}^T \1{I_t=i, \Nplay[][i]{t}{\gamma}< A}\\
 \leq  \sum_{t=1}^T \1{I_t=i, \Nplay[][i]{t-\tau:t}{1}< \gamma^{-\tau} A}
\eqsp, \label{eq:ducb:rhs1:prime}
\end{equation}
and apply the preceeding lemma with  $m= \gamma^{-\tau} A$.
\end{proof}
\subsection*{Acknowledgment}
The authors wish to thank Gilles Stoltz and Jean-Yves Audibert for stimulating discussions and for providing us with unpublished material.
\nocite{*}
\bibliography{./seqindivs}

\begin{thebibliography}{20}
\providecommand{\natexlab}[1]{#1}
\providecommand{\url}[1]{\texttt{#1}}
\expandafter\ifx\csname urlstyle\endcsname\relax
  \providecommand{\doi}[1]{doi: #1}\else
  \providecommand{\doi}{doi: \begingroup \urlstyle{rm}\Url}\fi

\bibitem[Agrawal(1995)]{agrawal:1995}
R.~Agrawal.
\newblock Sample mean based index policies with {$O(\log n)$} regret for the
  multi-armed bandit problem.
\newblock \emph{Adv. in Appl. Probab.}, 27\penalty0 (4):\penalty0 1054--1078,
  1995.
\newblock ISSN 0001-8678.

\bibitem[Audibert et~al.(2007)Audibert, Munos, and
  Szepesvári]{audibert:munos:szepesvari:2007}
J-Y. Audibert, R.~Munos, and A.~Szepesvári.
\newblock Tuning bandit algorithms in stochastic environments.
\newblock In \emph{Algorithmic Learning Theory}, 2007.

\bibitem[Auer et~al.(2002/03)Auer, Cesa-Bianchi, Freund, and
  Schapire]{auer:cesabianchi:freund:schapire:2002}
P.~Auer, N.~Cesa-Bianchi, Y.~Freund, and R.E. Schapire.
\newblock The nonstochastic multiarmed bandit problem.
\newblock \emph{SIAM J. Comput.}, 32\penalty0 (1):\penalty0 48--77
  (electronic), 2002/03.
\newblock ISSN 0097-5397.

\bibitem[Auer(2002)]{auer:2002}
Peter Auer.
\newblock Using confidence bounds for exploitation-exploration trade-offs.
\newblock \emph{J. Mach. Learn. Res.}, 3\penalty0 (Spec. Issue Comput. Learn.
  Theory):\penalty0 397--422, 2002.
\newblock ISSN 1532-4435.

\bibitem[Auer et~al.(2002)Auer, Cesa-Bianchi, and Fischer]{auer02finitetime}
Peter Auer, Nicol{\`{o}} Cesa-Bianchi, and Paul Fischer.
\newblock Finite-time analysis of the multiarmed bandit problem.
\newblock \emph{Machine Learning}, 47\penalty0 (2/3):\penalty0 235--256, 2002.
\newblock URL \url{citeseer.ist.psu.edu/auer00finitetime.html}.

\bibitem[Cesa-Bianchi and Lugosi(1999)]{cesabianchi:lugosi:1999}
Nicol{\`o} Cesa-Bianchi and G{\'a}bor Lugosi.
\newblock On prediction of individual sequences.
\newblock \emph{Ann. Statist.}, 27\penalty0 (6):\penalty0 1865--1895, 1999.
\newblock ISSN 0090-5364.

\bibitem[Cesa-Bianchi and Lugosi(2006)]{cesabianchi:lugosi:2006}
Nicolo Cesa-Bianchi and Gabor Lugosi.
\newblock \emph{Prediction, Learning, and Games}.
\newblock Cambridge University Press, New York, NY, USA, 2006.
\newblock ISBN 0521841089.

\bibitem[Cesa-Bianchi et~al.(2006)Cesa-Bianchi, Lugosi, and
  Stoltz]{cesabianchi:lugosi:stoltz:2006}
Nicol{\`o} Cesa-Bianchi, G{\'a}bor Lugosi, and Gilles Stoltz.
\newblock Regret minimization under partial monitoring.
\newblock \emph{Math. Oper. Res.}, 31\penalty0 (3):\penalty0 562--580, 2006.
\newblock ISSN 0364-765X.

\bibitem[Cesa-Bianchi et~al.(2008)Cesa-Bianchi, Lugosi, and
  Stoltz]{cesabianchi:lugosi:stoltz:2008}
Nicol{\`o} Cesa-Bianchi, G{\'a}bor Lugosi, and Gilles Stoltz.
\newblock Competing with typical compound actions, 2008.

\bibitem[Devroye et~al.(1996)Devroye, Gy{\"o}rfi, and Lugosi]{MR1383093}
L.~Devroye, L.~Gy{\"o}rfi, and G.~Lugosi.
\newblock \emph{A probabilistic theory of pattern recognition}, volume~31 of
  \emph{Applications of Mathematics (New York)}.
\newblock Springer-Verlag, New York, 1996.
\newblock ISBN 0-387-94618-7.

\bibitem[Freund and Schapire(1997)]{freund:schapire:robert:1995}
Yoav Freund and Robert~E. Schapire.
\newblock A decision-theoretic generalization of on-line learning and an
  application to boosting.
\newblock \emph{J. Comput. System Sci.}, 55\penalty0 (1, part 2):\penalty0
  119--139, 1997.
\newblock ISSN 0022-0000.
\newblock Second Annual European Conference on Computational Learning Theory
  (EuroCOLT '95) (Barcelona, 1995).

\bibitem[Hartland et~al.(2006)Hartland, Gelly, Baskiotis, Teytaud, and
  Sebag]{hartland:gelly:baskiotis:teytaud:sebag:2006}
C.~Hartland, S.~Gelly, N.~Baskiotis, O.~Teytaud, and M.~Sebag.
\newblock Multi-armed bandit, dynamic environments and meta-bandits, 2006.
\newblock NIPS-2006 workshop, Online trading between exploration and
  exploitation, Whistler, Canada.

\bibitem[Herbster and Warmuth(1998)]{herbster:warmuth:1998}
M.~Herbster and M.K. Warmuth.
\newblock Tracking the best expert.
\newblock \emph{Machine Learning}, 32\penalty0 (2):\penalty0 151--178, 1998.
\newblock ISSN 0885-6125.
\newblock \doi{http://dx.doi.org/10.1023/A:1007424614876}.

\bibitem[Kocsis and Szepesv\'ari(2006)]{kocsis:szepesvari:2006}
L.~Kocsis and C.~Szepesv\'ari.
\newblock Discounted {UCB}.
\newblock 2nd PASCAL Challenges Workshop, Venice, Italy, April 2006.

\bibitem[Koulouriotis and Xanthopoulos(2008)]{koulouriotis:xanthopoulos:2008}
D.~E. Koulouriotis and A.~Xanthopoulos.
\newblock Reinforcement learning and evolutionary algorithms for non-stationary
  multi-armed bandit problems.
\newblock \emph{Applied Mathematics and Computation}, 196\penalty0
  (2):\penalty0 913--922, 2008.

\bibitem[Kulkarni and Lugosi(2000)]{kulkarni:lugosi:2000}
S.~R. Kulkarni and G.~Lugosi.
\newblock Finite-time lower bounds for the two-armed bandit problem.
\newblock \emph{IEEE Trans. Automat. Control}, 45\penalty0 (4):\penalty0
  711--714, 2000.
\newblock ISSN 0018-9286.

\bibitem[Lai et~al.(2007)Lai, El~Gamal, Jiang, and Poor]{lai:2007}
L.~Lai, H.~El~Gamal, H.~Jiang, and H.~V. Poor.
\newblock Cognitive medium access: Exploration, exploitation and competition,
  2007.
\newblock URL
  \url{http://www.citebase.org/abstract?id=oai:arXiv.org:0710.1385}.

\bibitem[Lai and Robbins(1985)]{lai:robbins:1985}
T.~L. Lai and Herbert Robbins.
\newblock Asymptotically efficient adaptive allocation rules.
\newblock \emph{Adv. in Appl. Math.}, 6\penalty0 (1):\penalty0 4--22, 1985.
\newblock ISSN 0196-8858.

\bibitem[Slivkins and Upfal(2008)]{slivkins:upfal:2008}
A.~Slivkins and E.~Upfal.
\newblock Adapting to a changing environment: the brownian restless bandits,
  2008.
\newblock URL \url{http://research.microsoft.com/users/slivkins/}.
\newblock submitted to COLT 2008.

\bibitem[Whittle(1988)]{whittle:1988}
P.~Whittle.
\newblock Restless bandits: activity allocation in a changing world.
\newblock \emph{J. Appl. Probab.}, Special Vol. 25A:\penalty0 287--298, 1988.
\newblock ISSN 0021-9002.
\newblock A celebration of applied probability.

\end{thebibliography}

\clearpage
\end{document}